\documentclass[11 pt,reqno]{amsart}
\usepackage{amssymb}
\usepackage{amsmath}
\usepackage{amscd}

\input{diagrams}

\newtheorem {theorem}{Theorem}[section]
\newtheorem {proposition}[theorem]{Proposition}
\newtheorem {corollary}[theorem]{Corollary}
\newtheorem {lemma}[theorem]{Lemma}
\newtheorem {definition}[theorem]{Definition}

\numberwithin {equation}{section}
\renewcommand {\proof}{{\sc Proof.}\ }

\newcommand {\remark}
{\addtocounter{theorem}{1}{\bf Remark \thetheorem.}\ }

\newcommand {\halmos}{$\blacksquare$}

\newcommand {\IC}{\mathbb{C}}
\newcommand {\IN}{\mathbb{N}}                          
\newcommand {\IZ}{\mathbb{Z}}

\renewcommand {\c}{\mathfrak c}
\newcommand {\g}{\mathfrak{g}}
\newcommand {\h}{\mathfrak h}
\renewcommand {\ll}{\mathfrak l}
\newcommand {\n}{\mathfrak n}
\newcommand {\p}{\mathfrak p}

\renewcommand {\SS}{\mathfrak S}

\newcommand {\ad}{\operatorname{ad}}

\newcommand {\Aut}{\operatorname{Aut}}
\newcommand {\End}{\operatorname{End}}

\newcommand {\id}{\operatorname{id}}

\newcommand {\extg}[1]{\bigwedge^{#1}\g}

\newcommand {\fml}{[\negthinspace[\hbar]\negthinspace]}

\newcommand {\gl}{\mathfrak{gl}}
\newcommand {\Ug}{U\g}
\newcommand {\UgD}{U\g_D}

\renewcommand {\>}{\rangle}

\newcommand {\dH}{d_H}

\newcommand {\ie}{{\it i.e., }}
\newcommand {\eg}{{\it e.g.}, }
\newcommand {\fd}{finite--dimensional }
\newcommand {\lhs}{left--hand side }
\newcommand {\rhs}{right--hand side }
\renewcommand {\ss}{semi--simple }
\newcommand {\wrt}{with respect to }

\newcommand {\ol}{\overline}
\newcommand {\ul}{\underline}
\newcommand {\wt}{\widetilde}
\newcommand {\wh}{\widehat}

\newcommand {\CE}{Chevalley--Eilenberg }
\newcommand {\HC}{Harish--Chandra }
\newcommand {\Ho}{Hochschild }

\newcommand {\veps}{\varepsilon}

\newcommand {\ok}{^{\otimes k}}

\newcommand {\aand}{\qquad\text{and}\qquad}

\newcommand {\half}[1]{\frac{#1}{2}}

\newcommand {\cow}[1]{\lambda^{\vee}_{#1}}

\newcommand {\bas}{bialgebras }
\newcommand {\qtqba}{quasitriangular quasibialgebra }
\newcommand {\qcqtqba}{quasi--Coxeter quasitriangular quasibialgebra }
\newcommand {\qcqtqbas}{quasi--Coxeter quasitriangular quasibialgebras }

\newcommand {\YB}{\operatorname{YB}}

\newcommand {\Bg}{B_\g}
\newcommand {\Dg}{D_\g}

\newcommand {\schl}{[\negthinspace[}
\newcommand {\schr}{]\negthinspace]}

\newcommand {\cD}{\c_{D}}
\newcommand {\gD}{\g_{D}}

\newcommand {\lD}{\ll_{D}}
\newcommand {\nD}{\n_{D}}
\newcommand {\onD}{\ol{\n}_{D}^{-}}

\newcommand {\rg}{r_{\g}}
\newcommand {\rD}{r_{\gD}}

\newcommand {\Alt}{\operatorname{Alt}}
\newcommand {\Pent}{\operatorname{Pent}}

\begin{document}

\title{Cohomological construction of relative twists}
\author[V. Toledano Laredo]{Valerio Toledano Laredo}
\address{
Universit\'e Pierre et Marie Curie--Paris 6\\
Institut de Math\'ematiques de Jussieu, UMR 7586\\
Case 191\\
16 rue Clisson\\
F--75013 Paris}
\date{June 2005}
\email{toledano@math.jussieu.fr}
\begin{abstract}
Let $\g$ be a complex, semi--simple Lie algebra, $\h\subset\g$ a Cartan
subalgebra and $D$ a subdiagram of the Dynkin diagram of $\g$. Let
$\gD\subset\lD\subseteq\g$ be the corresponding semi--simple and Levi
subalgebras and consider two invariant solutions $\Phi\in(\Ug^{\otimes
3}\fml)^\g$ and $\Phi_D\in(\UgD^{\otimes 3}\fml)^{\gD}$ of the pentagon
equation for $\g$ and $\gD$ respectively. Motivated by the theory of
\qcqtqbas \cite{TL3}, we study in this paper the existence of a {\it relative
twist}, that is an element $F\in(\Ug^{\otimes 2}\fml)^{\lD}$ such that the
twist of $\Phi$ by $F$ is $\Phi_D$. Adapting the method of Donin and
Shnider \cite{DS}, who treated the case of an empty $D$, so that $\lD
=\h$ and $\Phi_D=1^{\otimes 3}$, we give a cohomological construction
of such an $F$ under the assumption that $\Phi_D$ is the image of $\Phi$
under the generalised \HC homomorphism
$(\Ug^{\otimes 3})^{\lD}\rightarrow(\UgD^{\otimes 3})^{\gD}$.
We also show that $F$ is unique up to a gauge transformation if $\lD$
is of corank 1 or $F$ satisfies $F^\Theta=F^{21}$ where $\Theta\in
\Aut(\g)$ is an involution acting as $-1$ on $\h$.
\end{abstract}
\maketitle

\section{Introduction}

Let $\g$ be a complex, semi--simple Lie algebra, $\h\subset\g$ a
Cartan subalgebra, $R_\g=\{\alpha\}\subset\h^*$ the corresponding
root system and $\Dg$ the Dynkin diagram of $\g$ relative to a
choice $\alpha_1,\ldots,\alpha_n\in\h^*$ of simple roots of $\g$.
Let $D\subseteq\Dg$ be a subdiagram of the Dynkin diagram of
$\g$ and denote by
$$\gD\subseteq\lD\subseteq\g$$
the corresponding diagrammatic subalgebra, \ie the \ss
subalgebra generated by the root vectors corresponding
to the simple roots in $D$, and Levi subalgebra $\lD=
\gD+\h$ respectively. Note that
$$\lD=\gD\oplus\cD$$
where the centre $\cD$ of $\lD$ is spanned by the fundamental
coweights $\cow{j}\in\h$, with $j$ such that $\alpha_j\notin D$.\\

Let $\hbar$ be a formal variable and consider two fixed, invariant
elements 
$$
\Phi\in 1+\hbar^2\left(\Ug^{\otimes 3}\fml\right)^{\g}
\quad\text{and}\quad
\Phi_D\in 1+\hbar^2\left(\UgD^{\otimes 3}\fml\right)^{\gD}
$$
satisfying the pentagon equation
\begin{equation}\label{eq:pentagon}
\id^{\otimes 2}\otimes\Delta(\Psi)
\cdot
\Delta\otimes\id^{\otimes 2}(\Psi)
=
1\otimes\Psi\cdot
\id\otimes\Delta\otimes\id(\Psi)\cdot
\Psi\otimes 1
\end{equation}

We shall be concerned in this paper with the cohomological solution
of the following {\it relative twist equation}
\begin{equation}\label{eq:twist}
(\Phi)_F
:=
1\otimes F\cdot\id\otimes\Delta(F)\cdot
\Phi
\cdot\Delta\otimes\id(F^{-1})\cdot F^{-1}\otimes 1
=\Phi_D
\end{equation}
\wrt an element $F$ which is invariant under the adjoint action of $\lD$
$$F\in 1+\hbar\left(\Ug^{\otimes 2}\fml\right)^{\lD}$$

Our motivation for studying \eqref{eq:twist} comes from the theory
of \qcqtqbas \cite{TL3}. These are, informally speaking, \bas which
carry representations of both Artin's braid groups $B_n$ and the
generalised braid group $\Bg$ of type $\g$ on the tensor products
of their \fd modules. One of the main results in \cite{TL3} is the
rigidity of \qcqtqba structures on $\Ug\fml$. In conjunction with
the results of \cite{TL4}, this shows in particular that the monodromy
of the Casimir connection introduced in \cite{MTL} is described by
Lusztig's quantum Weyl group operators \cite{Lu}, thus proving a
conjecture formulated independently by the author \cite{TL1,TL2}
and De Concini (unpublished). The rigidity result of \cite{TL3} depends
on Drinfeld's uniqueness theorem for \qtqba deformations of $\Ug$
\cite{Dr2} and on the uniqueness, up to gauge transformations, of
solutions of \eqref{eq:twist} when $\lD$ is of corank 1.\\

Rather than incorporating the required uniqueness result into \cite{TL3},
we decided to study the existence of solutions of \eqref{eq:twist} as well
and present our results in a separate publication. These may in fact be 
of independent interest since the relevant deformation complex turns
out to be a perturbation of the \CE complex for a suitable, non--coboundary
Lie algebra structure on $\g^*$. Our method is very close to that of
Donin--Shnider \cite[\S 3]{DS} who solved the equation \eqref{eq:twist}
when $D$ is empty, so that $\lD=\h$ and $\Phi_D=1^{\otimes 3}$, and
$\Phi$ satisfies in addition
\begin{equation}\label{eq:intro constraints}
\Phi^{321}=\Phi^{-1}
\quad\text{and}\quad
\Phi^{\Theta}=\Phi
\end{equation}
where $\Theta\in\Aut(\g)$ is an involution acting as $-1$ on $\h$. The possibility
of laddering down, that is solving \eqref{eq:twist} only when $|\Dg\setminus D|=1$
allows us to bypass the use of \eqref{eq:intro constraints} and to construct in
\S\ref{ss:coh existence} a suitable $F$ under the sole assumption that $\Phi
_D$ is the projection of $\Phi$ \wrt the generalised \HC homomorphism
$(\Ug^{\otimes 3})^{\lD}\rightarrow(\UgD^{\otimes 3})^{\gD}$ defined in
\S \ref{ss:HC}. Our proof proceeds along the lines of Donin and Shnider's,
the main difference being in the cohomology theory needed to deal with
secondary obstructions, which is defined and computed in \S\ref{ss:CE}.
The uniqueness of solutions of \eqref{eq:twist} is obtained in
\S\ref{ss:coh uniqueness} under the weaker assumption that the infinitesimal
of $\Phi$ projects onto that of $\Phi_D$ and that either $\lD$ is of corank 1
or $F$ satisfies $F^\Theta=F^{21}$. Section \ref{ss:standard YB} contains
some standard material on the classical Yang--Baxter equations.\\

\remark A non--cohomological proof of the existence of $F$ may be
given in the case where $\Phi$ and $\Phi_D$ are Lie associators
by adapting Etingof and Kazhdan's method \cite{EK}. The latter
corresponds to the case when $\lD=\h$ but can be modified by
replacing the Verma modules used in \cite{EK} by their generalised
counterparts obtained by inducing from the parabolic subalgebra
$\p_D\subset\g$ corresponding to $D$.

\section{Generalised \HC homomorphisms}
\label{ss:HC}

For each $k\geq 1$, we define below an algebra homomorphism
$$\pi_D^k:(\Ug\ok)^{\cD}\longrightarrow U\lD\ok$$
which restricts to the identity on $U\lD\ok$ and is equivariant \wrt
adjoint action of $\lD$. For $D=\emptyset$ and $k=1$, $\pi_D^k$
is the \HC homomorphism $\pi:\Ug^\h\rightarrow U\h$. The
definition of $\pi_D^k$ is similar to that of $\pi$, see \eg \cite
[\S 7.4.1--7.4.3]{Di} which we follow closely. Write
$$\g=\nD^{-}\oplus\lD\oplus\nD^{+}$$
where the nilpotent subalgebras $\nD^{\pm}$ are spanned by
the roots vectors $e_{\alpha},f_{\alpha}$ respectively, with
$\alpha$ ranging over the positive roots of $\g$ not lying in
the root system $R_D$ of $\gD$. Set
$$I_k=
(\Ug\ok)^{\cD}\bigcap\medspace
\sum_{i=1}^k\Ug\ok\cdot(\nD^{+})_i$$
where, for $y\in \Ug$,
$$y_i=1^{\otimes(i-1)}\otimes y\otimes 1^{\otimes(k-i)}\in\Ug\ok$$

\begin{proposition}\label{pr:HC}\hfill
\begin{enumerate}
\item $I_k=(\Ug\ok)^{\cD}\bigcap\medspace
\sum_{i=1}^k(\nD^{-})_i\cdot\Ug\ok$.
\item $I_k$ is a two--sided ideal in $(\Ug\ok)
^{\cD}$ invariant under the adjoint action of $\lD$.
\item $(\Ug\ok)^{\cD}=I_k\oplus U\lD\ok$.
\end{enumerate}
\end{proposition}
\proof (i) By the PBW theorem,
$$\Ug\ok\cong
{U\nD^{-}}\ok\otimes
U\lD\ok\otimes
{U\nD^{+}}\ok$$
is spanned by the monomials
$$
u(q_i^j;x;p_i^j)=
f_{\beta_1,1}^{q_1^1}\cdots f_{\beta_{m},1}^{q_{m}^1}
\cdots
f_{\beta_1,k}^{q_1^k}\cdots f_{\beta_{m},k}^{q_{m}^k}
\cdot x \cdot
e_{\beta_1,1}^{q_1^1}\cdots e_{\beta_{m},1}^{q_{m}^1}
\cdots
e_{\beta_1,k}^{q_1^k}\cdots e_{\beta_{m},k}^{q_{m}^k}
$$
where $x\in U\lD\ok$, $\beta_1,\ldots,\beta_{m}$ are the
positive roots in $R_{\g}\setminus R_D$ and $q_i^j,p_i
^j\in\IN$. Let $\imath^*:\h^{*}\rightarrow\cD^{*}$ be the restriction
map. Since $u(q_i^j;x;p_i^j)$ has weight $\imath^*\sum_
{i,j}p_{i,j}\beta_j-\imath^*\sum_{i,j}q_{i,j}\beta_j$ for the
adjoint action of $\cD$, $(\Ug\ok)^{\cD}$ is spanned by
the $u(q_i^j;x;p_i^j)$ such that 
$$
\imath^*\sum_{i,j}p_{i,j}\beta_j=
\imath^*\sum_{i,j}q_{i,j}\beta_j
$$
Note that, since each $\beta_j$ restricts on $\cD$ to a
non--trivial linear combination of the simple roots
$\alpha_\ell\notin D$, with non--negative coefficients, $
\imath^*\sum_{i,j}p_{i,j}\beta_j=0$ iff $\sum_{i,j}p_{i,j}=0$.
It follows that
\begin{equation*}
\begin{split}
I_k
&=
\<u(q_i^j;x;p_i^j)\in(\Ug\ok)^{\cD}\>
_{\sum_{i,j}p_i^j>0}\\
&=
\<u(q_i^j;x;p_i^j)\in(\Ug\ok)^{\cD}\>
_{\sum_{i,j}q_i^j>0}\\
&=
(\Ug\ok)^{\cD}\bigcap\medspace
\sum_{i=1}^k(\nD^{-})_i\cdot\Ug\ok
\end{split}
\end{equation*}
as claimed. (ii) $I_k$ is a left ideal by definition and,
by (i), it is also a right ideal. It is moreover invariant
under the adjoint action of $\lD$ since $\nD^{\pm}$
are. (iii) is now obvious \halmos 

\begin{corollary}
The projection $\pi_D^k$ of $(\Ug\ok)^{\cD}$ onto $U\lD\ok$
defined by the ideal $I_k$ is equivariant for the adjoint action
of $\lD$ and therefore gives rise to the following commutative
diagram of algebra homomorphisms
$$
\begin{diagram}[labelstyle=\scriptstyle,balance]
\left(\Ug\ok\right)^{\cD} & \rTo^{\pi^k_D} & U\lD\ok		     
& 
\rTo & \UgD\ok\\
\bigcup				   &			& \bigcup	   		     
&
	    & \bigcup	    \\
\left(\Ug\ok\right)^{\lD} & \rTo^{\pi^k_D} & 
\left(U\lD\ok\right)^{\lD}& \rTo & \left(\UgD\ok\right)^{\gD}
\end{diagram}
$$
where the rightmost horizontal arrows are induced by the Lie
algebra projection $\lD\rightarrow\gD$.
\end{corollary}

\begin{definition}\label{de:HC}
We denote the composition of the horizontal arrows by $\ol{\pi}
_D^k$ and refer to $\ol{\pi}_D^k$ or $\pi_D^k$ as
generalised \HC homomorphisms.
\end{definition}

Note that $\pi_D^k$ and $\ol{\pi}_D^k$ are equivariant
under the natural action of the symmetric group $\SS_k$.
We record for later use the following two results

\begin{proposition}\label{pr:HC face}
For any $i,l\leq k$, $x\in(\Ug\ok)^{\cD}$ and $y\in(\Ug^{\otimes l})
^{\cD}$, one has
\begin{gather}
\id^{\otimes i}\otimes\Delta\otimes\id^{\otimes (k-i-1)}
\circ
\pi_D^k(x)
=
\pi_D^{k+1}
\circ
\id^{\otimes i}\otimes\Delta\otimes\id^{\otimes (k-i-1)}(x)
\label{eq:face 1}
\\
1^{\otimes i}\otimes\pi_D^{l}(y)\otimes 1^{\otimes(k-l-i)}
=
\pi_D^k
(1^{\otimes i}\otimes y\otimes 1^{\otimes(k-i-l)})
\label{eq:face 2}
\end{gather}
These identities remain valid if $\pi_D^k,\pi_D^{k+1}$
and $\pi_D^{l}$ are replaced by $\ol{\pi}_D^k,\ol{\pi}
_D^{k+1}$ and $\ol{\pi}_D^{l}$ respectively.
\end{proposition}
\proof Let 
$$
\gamma=
\id^{\otimes i}\otimes\Delta\otimes\id^{\otimes (k-i-1)}:
\Ug\ok\rightarrow\Ug^{\otimes(k+1)}
$$
Since $\gamma$ is equivariant for the adjoint action of $\g$,
it maps $(\Ug\ok)^{\cD}$ to $(\Ug^{\otimes (k+1)})^{\cD}$ so
that the \rhs of \eqref{eq:face 1} is well--defined.
One readily checks that
$$
\gamma(I_k)\subset I_{k+1}
\quad\text{and that}\quad
\gamma(U\lD\ok)\subset U\lD^{\otimes(k+1)}
$$
so that \eqref{eq:face 1} holds. \eqref{eq:face 2} is proved
in the same way. The fact that these identities hold when 
$\pi_D^j$ is replaced by $\ol{\pi}_D^j$ throughout
follows from the fact that $\ol{\pi}_D^j=\pi^{\otimes 
j}\circ\pi_D^j$ where $\pi:U\lD\rightarrow\Ug_D$ is
a Hopf algebra homomorphism \halmos

\begin{corollary}\label{co:HC Ho}
Let $\dH:\Ug^{\otimes k}\rightarrow\Ug^{\otimes(k+1)}$ be
the \Ho differential given by
\begin{equation}\label{eq:Hochschild diff}
\dH\medspace x
=1\otimes x+
\sum_{i=1}^k(-1)^i
\id^{\otimes(i-1)}\otimes\Delta\otimes\id^{\otimes(k-i)}
\medspace(x)
+(-1)^{k+1} x\otimes 1
\end{equation}
Then,
$$
\dH\circ\pi_D^k=\pi_D^{k+1}\circ\dH
\quad\text{and}\quad
\dH\circ\ol{\pi}_D^k=\ol{\pi}_D^{k+1}\circ\dH
$$
\end{corollary}

\section{Classical Yang--Baxter equations}
\label{ss:standard YB}

We review below some well--known results on the classical
Yang--Baxter equations due to Drinfeld \cite{Dr1}.\\

Define the classical Yang--Baxter map $\YB:\g^{\otimes 2}
\otimes\g^{\otimes 2}\rightarrow\g^{\otimes 3}$ by
$$
\YB(r,s)=
[r^{12},s^{13}+s^{23}]+[r^{13},s^{23}]
+
[s^{12},r^{13}+r^{23}]+[s^{13},r^{23}]
$$
Identify the exterior algebra $\bigwedge\g$ with its image
in the tensor algebra $T\g$ via the antisymmetrisation map
$$
X_1\wedge\ldots\wedge X_k
\longrightarrow
\Alt_k(X_1\otimes\cdots\otimes X_k)=
\frac{1}{k!}\sum_{\sigma\in\SS_k}(-1)^{\sigma(1)}
X_{\sigma(1)}\otimes\cdots\otimes X_{\sigma(k)}
$$
One readily checks that if $r,s\in\bigwedge^2\g$, then
\begin{equation}\label{eq:Alt=YB}
\YB(r,s)=6\Alt_3[r^{12},s^{13}]
\end{equation}

\begin{lemma}\label{le:YB=sch}
If $r=r_1\wedge r_2,s=s_1\wedge s_2\in\bigwedge
^2\g$, then
\begin{equation}\label{eq:YB=sch}
\YB(r,s)=
\frac{3}{2}
\sum_{1\leq i,j\leq 2}
[r_i,s_j]\wedge r_{3-i}\wedge s_{3-j}
\end{equation}
\end{lemma}
\proof We have
\begin{equation*}
\begin{split}
4[r^{12},s^{13}]
&=
 [r_1,s_1]\otimes r_2\otimes s_2
-[r_1,s_2]\otimes r_2\otimes s_1\\
&\phantom{=}
-[r_2,s_1]\otimes r_1\otimes s_2
+[r_2,s_2]\otimes r_1\otimes s_1
\end{split}
\end{equation*}
Antisymmetrising both sides and using \eqref{eq:Alt=YB},
we find \eqref{eq:YB=sch} \halmos\\

Let $(\cdot,\cdot)$ be a non--degenerate, ad--invariant,
symmetric bilinear form on $\g$ and let $\Omega=\sum_i
x_i\otimes x^i$, where $\{x_i\},\{x^i\}$ are dual basis of
$\g$ with respect to $(\cdot,\cdot)$, be the corresponding
symmetric, invariant tensor in $\g\otimes\g$. It is well--know
that $[\Omega_{12},\Omega_{23}]$ lies in $(\bigwedge^3\g)
^\g$ and generates it if $\g$ is simple. Let
\begin{equation}\label{eq:standard r}
\rg=
\sum_{\alpha\succ 0}\half{(\alpha,\alpha)}\cdot
e_{\alpha}\wedge f_{\alpha}\in\bigwedge^2\g
\end{equation}
where $e_{\alpha}\in\g_{\alpha}$, $f_{\alpha}\in\g_{-\alpha}$
are root vectors such that $[e_{\alpha},f_{\alpha}]=h_{\alpha}$
so that
\begin{equation}\label{eq:dual}
(e_{\alpha},f_{\alpha})=
\half{1}([h_{\alpha},e_{\alpha}],f_{\alpha})=
\half{1}(h_{\alpha},[e_{\alpha},f_{\alpha}])=
\half{1}(h_{\alpha},h_{\alpha})=
\frac{2}{(\alpha,\alpha)}
\end{equation}
By the following result, $\rg$ is a solution of the modified
classical Yang--Baxter equation (MCYBE), that is the equation
\begin{equation}\label{eq:MCYBE}
[\rg^{12},\rg^{23}+\rg^{13}]+[\rg^{13},\rg^{23}]
\in(\bigwedge^3\g)^{\g}
\end{equation}

\begin{proposition}[Drinfeld]\label{pr:MCYBE}
$$\YB(\rg,\rg)
=\frac{1}{2}
[\Omega_{12},\Omega_{23}]$$
\end{proposition}

\remark We shall refer to $\rg$ given by \eqref{eq:standard r}
as the standard (Drinfeld) solution of the MCYBE corresponding
to the bilinear form $(\cdot,\cdot)$.

\section{Classical $r$--matrices and \CE cohomology}
\label{ss:CE}

\subsection{}

The aim of this section is to compute the cohomology
of the complex
$$
((\bigwedge\g)^{\gD},d)
\quad\text{where}\quad
d=\schl\rg-\rD,\cdot\schr
$$
is given by the Schouten bracket with the difference of the
standard solutions of the modified classical Yang--Baxter
equations for $\g$ and $\gD$ respectively.\\

The computation is carried out by identifying $d$ with a 
perturbation of the \CE differential on $\bigwedge\g =
\bigwedge(\g^{*})^{*}$ induced by a suitable Lie algebra
structure on $\g^{*}$. When $D=\emptyset$, so that $\lD=\h$,
this identification is well--known and follows readily from the
fact that the relevant Lie algebra structure on $\g^*$ is given
in terms of the cobracket $\delta:\g\rightarrow\g\wedge\g$
defined by
$$\delta(X)=\schl\rg,X\schr=-\ad(X)\rg$$
When $D\neq\emptyset$ the relevant Lie algebra structure
on $\g^*$ is described in \S \ref{sss:g*} and is not of
coboundary type. We begin with a few reminders.

\subsection{}\label{sss:schouten}

Recall that the Schouten bracket
$$\schl\cdot,\cdot\schr :
\bigwedge^k\g\otimes\bigwedge^{l}\g
\rightarrow
\bigwedge^{k+l-1}\g$$
on the exterior algebra $\bigwedge\g$ is defined by
\begin{multline}\label{eq:schouten}
\schl
X_1\wedge\cdots\wedge X_k,
Y_1\wedge\cdots\wedge Y_{l}
\schr\\
=
\sum_{i,j}(-1)^{i+j}
[X_i,Y_j]\wedge
X_1\wedge\cdots\wedge\wh{X_i}\wedge\cdots\wedge X_k
\wedge
Y_1\wedge\cdots\wedge\wh{Y_j}\wedge\cdots\wedge Y_{l}
\end{multline}
The Schouten bracket satisfies
$$\schl\ul{X},\ul{Y}\schr=-(-1)^{(k-1)(l-1)}\schl\ul{Y},\ul{X}\schr$$
for any $\ul{X}\in\bigwedge^k\g$ and $\ul{Y}\in\bigwedge^{l}\g$
and
$$\schl\ul{X},\schl\ul{Y},\ul{Z}\schr\schr
=
\schl\schl\ul{X},\ul{Y}\schr,\ul{Z}\schr
+(-1)^{(k-1)(l-1)}
\schl\ul{Y},\schl\ul{X},\ul{Z}\schr\schr$$
for any such $\ul{X},\ul{Y}$ and $\ul{Z}\in\bigwedge\g$, and therefore
endows $\bigwedge\g$ with the structure of a $\IZ$--graded Lie algebra,
provided its grading is defined by
$$\deg(\bigwedge^k\g)=k-1$$

Moreover, since
$$\schl\ul{X},\ul{Y}\wedge\ul{Z}\schr=
\schl\ul{X},\ul{Y}\schr\wedge\ul{Z}+
(-1)^{(k-1)l}\ul{Y}\wedge\schl\ul{X},\ul{Z}\schr$$
for any  $\ul{X}\in\bigwedge^k\g$ and $\ul{Y}\in\bigwedge^l\g$,
the map $\ul{X}\rightarrow\schl\ul{X},\cdot\schr$ is a homomorphism
of $\bigwedge\g$ into the $\IZ$--graded Lie algebra of derivations
of the exterior algebra $\bigwedge\g$ endowed with its standard
grading.\\

Note that any $r\in\bigwedge^2\g$ defines a degree 1 derivation
$d_{r}=\schl r,\cdot\schr$ of $\bigwedge\g$. Its square is readily
computed from
$$d_{r}^2(\ul{Y})
=
\schl r,\schl r,\ul{Y}\schr\schr
=
\schl\schl r,r\schr,\ul{Y}\schr-
\schl r,\schl r,\ul{Y}\schr\schr
=
\schl\schl r,r\schr,\ul{Y}\schr-
d_{r}^2(\ul{Y})$$
Since $d_{r}^2$ is also an algebra derivation of $\bigwedge
\g$, and $\schl\ul{X},Y\schr=-\ad(Y)\ul{X}$ for any $\ul{X}
\in\bigwedge\g$ and $Y\in\g$, $d_{r}$ is a differential if,
and only if 
$$\schl r,r\schr\in(\bigwedge^3\g)^{\g}$$
and therefore, by lemma \ref{le:YB=sch}, iff $r$ is a solution
of the MCYBE \eqref{eq:MCYBE}.

\subsection{}\label{sss:rg-rD}

Let now
$$\rg\in\bigwedge^2\g
\qquad\text{and}\qquad
\rD\in\bigwedge^2\gD$$
be solutions of the MCYBE for $\g$ and $\gD$ respectively such
that $\rg-\rD$ is invariant under $\gD$. This is the case for example
if both $\rg$ and $\rD$ are the standard solutions \eqref{eq:standard r}
of MCYBE relative to a non--degenerate, ad--invariant bilinear form
$(\cdot,\cdot)$ on $\g$ and its restriction to $\gD$ respectively. Indeed,
$\nD^{+}$ and $\nD^{-}$ are invariant under the adjoint action of
$\gD$ and $(\cdot,\cdot)$ yields a $\gD$--equivariant identification
$(\nD^{+})^{*}\cong \nD^{-}$ with respect to which
\begin{equation}\label{eq:rg-rD}
\rg-\rD=
\sum_{\alpha\in R_{\g}^{+}\setminus R_D}
\half{(\alpha,\alpha)}\cdot e_{\alpha}\wedge f_{\alpha}
\end{equation}
is the image in $\bigwedge^2(\nD^{+}\oplus\nD^{-})$ of
$$\id_{\nD^{+}}\in
\End(\nD^{+})\cong
\nD^{+}\otimes\nD^{-}\subset
(\nD^{+}\oplus\nD^{-})^{\otimes 2}$$
under the projection $(\nD^{+}\oplus\nD^{-})^{\otimes 2}
\rightarrow\bigwedge^2(\nD^{+}\oplus\nD^{-})$. We shall
need the following simple

\begin{lemma}\label{le:simple}
For any $\ul{X}=X_1\wedge\cdots\wedge X_k\in\bigwedge^k
\g$, the following holds on $\bigwedge\g$
\begin{equation}\label{eq:simple}
\schl\ul{X},\cdot\schr
=
(-1)^{k-1}\cdot\sum_{i=1}^k(-1)^{i-1}
e(X_1\wedge\cdots\wedge\wh{X_i}\wedge\cdots\wedge X_k)
\cdot
\ad(X_i)
\end{equation}
where $e(\ul{Y})$ is exterior multiplication by $\ul{Y}$. In
particular, if $\ul{X}\in\bigwedge\gD$ and $\ul{Y}\in(\bigwedge
\g)^{\gD}$, then $\schl\ul{X},\ul{Y}\schr=0$.
\end{lemma}

\begin{proposition}\label{pr:differential}
Let $\rg,\rD$ be solutions of the MCYBE for $\g,\gD$ respectively 
such that $\rg-\rD$ is invariant under $\gD$. Then
\begin{enumerate}
\item $\schl\rg-\rD,\cdot\schr$ leaves $(\bigwedge\g)^{\gD}$
invariant.
\item Its restriction to $(\bigwedge\g)^{\gD}$ coincides with
that of $\schl\rg,\cdot\schr$ and is therefore a differential.
\item $\schl\rg-\rD,\rg-\rD\schr=\schl\rg,\rg\schr-\schl\rD,
\rD\schr$.
\end{enumerate}
\end{proposition}
\proof (i) Since $\rg-\rD$ is invariant under $\gD$, and the
Schouten bracket is equivariant for the adjoint action of $\g$,
$\schl\rg-\rD,\cdot\schr$ leaves $(\bigwedge\g)^{\gD}$ invariant.
(ii) By lemma \ref{le:simple}, $\schl\rD,\ul{Y}\schr=0$ for any
$\ul{Y}\in(\bigwedge\g)^{\gD}$ so that
\begin{equation}\label{eq:res}
\schl\rg-\rD,\ul{Y}\schr=\schl\rg,\ul{Y}\schr
\end{equation}
for any such $\ul{Y}$. (iii) Since $\rg-\rD$ is invariant under
$\gD$, we find, by \eqref{eq:res}
\begin{multline*}
\schl\rg-\rD,\rg-\rD\schr
=
\schl\rg,\rg-\rD\schr\\
=
\schl\rg,\rg\schr-\schl\rg-\rD,\rD\schr-\schl\rD,\rD\schr
=\schl\rg,\rg\schr-\schl\rD,
\rD\schr
\end{multline*}
as claimed \halmos\\

\remark Note that the proof of (iii) only uses the $\gD
$--invariance of $\rg-\rD$. Thus, if $\rg\in\bigwedge^
{2}\g$ is a solution of the MCYBE and $\rD\in\bigwedge^
{2}\gD$ is such that $\rg-\rD$ is invariant under $\gD$,
then $\rD$ is a solution of the MCYBE for $\gD$. Moreover,
if $\rg$ is the standard solution of the MCYBE then so is
$\rD$. Indeed, if $\pi\in\End(\bigwedge^2\g)$ is the
projection onto $\gD$--invariants, then
\begin{multline*}
\rg-\rD=
\pi(\rg-\rD)\\
=
\pi(\rg-\ol{\pi}_D^2(\rg)+\ol{\pi}_D^2(\rg)-\rD)=
\rg-\ol{\pi}_D^2(\rg)
\end{multline*}
where the last equality follows from the $\gD$--invariance
of $\rg-\ol{\pi}_D^2(\rg)$ and the fact that $\pi(\ol
{\pi}_D^2(\rg)-\rD)\in(\bigwedge^2\gD)^{\gD}=\{0\}$.
Thus, $\rD=\ol{\pi}_D^2(\rg)$ is the standard solution
of the MCYBE for $\gD$.

\subsection{}\label{sss:g*}

Identify $\g^{*}$ and $\g$ as vector spaces by using the bilinear
form $(\cdot,\cdot)$, and endow $\g^{*}$ with the following Lie
algebra structure
\begin{equation}\label{eq:Lie str}
\g^{*}=(\nD^{+}\oplus\onD)\rtimes(\gD\oplus\cD)
\end{equation}
where $\onD$ is $\nD^{-}$ with the opposite bracket, $\gD$
acts on $\nD^{\pm}$ by the adjoint action and $\cD$ acts on
$\nD^{\pm}$ by $\pm 1/2$ times the adjoint action. Denoting
the corresponding bracket on $\g^*$ by $[\cdot,\cdot]^*$, we
therefore have
\begin{equation}\label{eq:new bracket}
\begin{split}
[x,y]^*
&=
[x_D,y_D]+[x_D,y_{+}+y_{-}]+[x_{+}+x_{-},y_D]\\
&+
[x_{+},y_{+}]-[x_{-},y_{-}]+
\half{1}[x_{0},y_{+}-y_{-}]+
\half{1}[x_{+}-x_{-},y_{0}]
\end{split}
\end{equation}
where, $z_D\in\gD$, $z_{\pm}\in\nD^{\pm}$ and $z_{0}\in
\cD$ are the components of $z\in\g^{*}$ corresponding to
the decomposition \eqref{eq:Lie str}. Thus, $\gD$ is a Lie
subalgebra of $\g^{*}$ and its coadjoint action on $(\g^{*})
^{*}=\g$ coincides with its adjoint action on $\g$.\\

Let now $\delta\in\End(\bigwedge\g)$ be the differential
obtained by regarding $\bigwedge\g$ as the \CE complex of
$\g^{*}$. The following result identifies $\schl\rg-\rD,
\cdot\schr$ with a perturbation of $\delta$.

\begin{theorem}\label{th:sch=CE}
If $\rg$ and $\rD$ are the standard solutions of the MCYBE
corresponding to $(\cdot,\cdot)$ and its restriction to $\gD$
respectively, the following holds on $\bigwedge\g$
\begin{equation}\label{eq:sch=CE}
\schl\rg-\rD,\cdot\schr=
2\delta+e(v_i)\cdot[\ad(v_i)\cdot(1+2P_{+})]^{\wedge}
\end{equation}
where $\{v_i\},\{v^i\}$ are basis of $\gD$ dual with
respect to $(\cdot,\cdot)$, $P_{+}:\g\rightarrow\nD^{+}$
is the projection corresponding to the decomposition \eqref
{eq:Lie str} and $T\in\gl{}(\g)\rightarrow T^{\wedge}\in\gl{}
(\bigwedge\g)$ is the Lie algebra homomorphism given by
$$
T^{\wedge}\medspace X_1\wedge\ldots\wedge X_k
=
\sum_{i=1}^k
X_1\wedge\ldots\wedge TX_i\wedge\ldots\wedge X_k
$$
\end{theorem}
\proof It is sufficient to check \eqref{eq:sch=CE} on elements
of $\g\subset\bigwedge\g$ since both sides are degree 1 algebra
derivations of $\bigwedge\g$. In turn, it is easier to check
that the transposes of both sides coincide as maps $\bigwedge
^2\g\rightarrow\g$. By definition, $\delta^{t}=[\cdot,\cdot]
^{*}$. Since $e(v)^{t}=\imath(v)$ where $\imath(v)$ is the
contraction operator defined by
$$
\imath(v)\medspace Y_1\wedge\cdots\wedge Y_{l}
=
\sum_{i=1}^{l}(-1)^{i-1}(v,Y_i)\medspace
Y_1\wedge\cdots\wedge\wh{Y_i}\wedge\cdots\wedge Y_{l}
$$
and $\ad(X)^{t}=-\ad(X)$ for any $X\in\g$, we find, using
\eqref{eq:rg-rD} and \eqref{eq:simple}
$$
\schl\rg-\rD,\cdot\schr^{t}
=
\sum_{\alpha\in R_{\g}^{+}\setminus R_D}
\frac{(\alpha,\alpha)}{2}
\left(
\ad(e_{\alpha})\imath(f_{\alpha})-
\ad(f_{\alpha})\imath(e_{\alpha})
\right)
$$
which, applied to $u\wedge v\in\bigwedge^2\g$ yields
\begin{equation*}
\begin{split}
&
\sum_{\alpha\in R_{\g}^{+}\setminus R_D}
\frac{(\alpha,\alpha)}{2}
\left(
(f_{\alpha},u)[e_{\alpha},v]-(f_{\alpha},v)[e_{\alpha},u]-
(e_{\alpha},u)[f_{\alpha},v]+(e_{\alpha},v)[f_{\alpha},u]
\right)\\
&
\phantom{\sum(\alpha,\alpha)}
=
[u_{+},v]+[u,v_{+}]-[u_{-},v]-[u,v_{-}]\\
&
\phantom{\sum(\alpha,\alpha)}
=
2[u_{+},v_{+}]-2[u_{-},v_{-}]+
 [u_{0},v_{+}-v_{-}]+[u_{+}-u_{-},v_{0}]\\
&
\phantom{\sum(\alpha,\alpha)=}
+
[u_{+},v_D]+[u_D,v_{+}]-[u_{-},v_D]-[u_D,v_{-}]
\end{split}
\end{equation*}
Comparing with \eqref{eq:new bracket}, we see that this is
equal to
\begin{equation*}
\begin{split}
 &
2[u,v]^{*}-2[u_D,v_D]
-[u_D,v_{+}]-[u_{+},v_D]
-3[u_D,v_{-}]-3[u_{-},v_D]\\
=&
2[u,v]^{*}-[u_D,v]-[u,v_D]-2[u_D,v_{-}]-2[u_{-},v_D]\\
=&
2[u,v]^{*}-(1+2P_{-})([u_D,v]+[u,v_D])
\end{split}
\end{equation*}
where $P_{-}$ is the projection onto $\nD^{-}$ which commutes
with the adjoint action of $\gD$. Noting that, for $x,y\in\g$,
one has
$$
\ad(v^i)\imath(v_i)x\wedge y=
(v_i,x)[v^i,y]-(v_i,y)[v^i,x]=
[x_D,y]+[x,y_D]
$$
we therefore find
$$
\schl\rg-\rD,\cdot\schr^{t}=
2[\cdot,\cdot]^{*}-(1+2P_{-})\ad(v_i)\imath(v^i)
$$
which yields \eqref{eq:sch=CE} since $P_{-}^{t}=P_{+}$ \halmos

\subsection{}\label{sss:cohomology}

Note that $(\bigwedge\g)^{\lD}$ is a subcomplex of $\left(
(\bigwedge\g)^{\gD},\schl\rg-\rD,\cdot\schr\right)$ since
$\rg-\rD$ is of weight zero. Note also that the restriction
of $\schl\rg-\rD,\cdot\schr$ to
$$
(\bigwedge\lD)^{\lD}=
(\bigwedge\gD)^{\gD}\wh{\otimes}\bigwedge\cD
\subset(\bigwedge\g)^{\lD}
$$
is zero since $\rg-\rD$ is invariant under $\lD$.

\begin{theorem}\label{th:sch coh}
If $\rg,\rD$ are the standard solutions of the MCYBE for $\g,
\gD$ respectively, the inclusions
$$
\left((\bigwedge\lD)^{\lD},0\right)
\longrightarrow
\left((\bigwedge\g)^{\lD},\schl\rg-\rD,\cdot\schr\right)
\longrightarrow
\left((\bigwedge\g)^{\gD},\schl\rg-\rD,\cdot\schr\right)
$$
are quasi--isomorphisms.
\end{theorem}
\proof Denote $\schl\rg-\rD,\cdot\schr$ by $d$. It is sufficient
to find an $\lD$--equivariant, diagonalisable operator $C\in\End
(\bigwedge\g)$ with kernel $\bigwedge\lD$ and an $\lD$--equivariant
homotopy $h\in\End(\bigwedge\g)$ such $dh+hd=C$. Noting that $\cD
\subset\g^{*}$ acts on $\bigwedge\g$ via the coadjoint action with
non--negative weights only so that the corresponding subspace of
invariants is precisely $\bigwedge\lD$, we see that a suitable $C$ 
is given by the Casimir operator
$$
C=\ad^{*}(t_i)\ad^{*}(t^i)
$$
where $\ad^{*}$ is the coadjoint action of $\g^{*}$ on $\bigwedge
\g$ and $\{t_i\},\{t^i\}$ are dual basis of $\cD$ with respect
to $(\cdot,\cdot)$. We claim that
$$
h=\ad^{*}(t_i)\imath(t^i)
$$
satisfies $dh+hd=2C$. It is well--known that $h$ satisfies $\delta
h+h\delta=C$, where $\delta\in\End(\bigwedge\g)$ is the \CE 
differential. Indeed,
\begin{equation*}
\begin{split}
\delta h+h\delta
&=
\delta\ad^{*}(t_i)\imath(t^i)+\ad^{*}(t_i)\imath(t^i)\delta\\
&=
\ad^{*}(t_i)(\delta\medspace\imath(t^i)+\imath(t^i)\delta)\\
&=
\ad^{*}(t_i)\ad^{*}(t^i)
\end{split}
\end{equation*}
where we have used the fact that $\delta$ is equivariant for
$\ad^{*}$ and the identity $\delta\imath(X)+\imath(X)\delta=
\ad^{*}(X)$, $X\in\g^{*}$. By theorem \ref{th:sch=CE}, it
therefore suffices to show that $h$ anticommutes with
$$
k=e(v_j)\cdot[(1+2P_{+})\cdot\ad(v^j)]^{\wedge}
$$
Bearing in mind the following identities for $X,Y\in\g$ and
$T\in\gl{}(\g)$
\begin{gather*}
\imath(X)e(Y)+e(Y)\imath(X)=(X,Y),\\
[T^{\wedge},\imath(X)]=\imath(TX)
\qquad\text{and}\qquad
[T^{\wedge},e(Y)]=e(TY)
\end{gather*}
and the fact that $\ad^{*}(X)^{\wedge}=\ad^{*}(X)$ for any 
$X\in\g^{*}$, we find
\begin{equation*}
\begin{split}
kh
&=
e(v_j)\cdot
[(1+2P_{+})\cdot\ad(v^j)]^{\wedge}\cdot
\ad^{*}(t_i)\imath(t^i)\\
&=
\ad^{*}(t_i)e(v_j)\cdot
[(1+2P_{+})\cdot\ad(v^j)]^{\wedge}\cdot
\imath(t^i)\\
&=
-\ad^{*}(t_i)\imath(t^i)
e(v_j)\cdot[(1+2P_{+})\cdot\ad(v^j)]^{\wedge}
\end{split}
\end{equation*}
as claimed \halmos\\

Since $(\bigwedge^i\gD)^{\gD}=0$ for $i=1,2$, we obtain
in particular the following 
\begin{corollary}\label{co:small cohomology}
\begin{align*}
H^1((\bigwedge\g)^{\gD};\schl\rg-\rD,\cdot\schr)
&\cong\cD \\
H^2((\bigwedge\g)^{\gD};\schl\rg-\rD,\cdot\schr)
&\cong\bigwedge^2\cD \\
H^3((\bigwedge\g)^{\gD};\schl\rg-\rD,\cdot\schr)
&\cong\bigwedge^3\cD\oplus(\bigwedge^3\gD)^{\gD}
\end{align*}
\end{corollary}

\section{Existence of twists}
\label{ss:coh existence}

\subsection{}

Let
$$\Phi\in 1^{\otimes 3}+\hbar^2\left(\Ug^{\otimes 3}\fml\right)^{\g}$$
be a solution of the pentagon equation \eqref{eq:pentagon}.
We shall need to assume that $\Phi$ is non--degenerate in
the sense defined below. Write
$$
\Phi=1^{\otimes 3}+\hbar^2\varphi\thickspace\mod\hbar^3
\quad\text{where}\quad
\varphi\in(\Ug^{\otimes 3})^{\g}
$$
Taking the coefficient of $\hbar^2$ in the pentagon
relation for $\Phi$, we find that $\dH\medspace\varphi=0$
where $\dH$ is the \Ho differential given by \eqref{eq:Hochschild diff}.
Thus,
\begin{equation}\label{eq:decomposition}
\Alt_3(\varphi)\in
(\bigwedge^3\g)^{\g}=
\bigoplus_i(\bigwedge^3\g_i)^{\g_i}
\end{equation}
where $\g_i$ are the simple factors of $\g$.

\begin{definition}\label{de:non deg}
$\Phi$ is a non--degenerate solution of the pentagon equation
if the components of $\Alt_3(\varphi)$ along the decomposition
\eqref{eq:decomposition} are all non--zero.
\end{definition}

Since each $(\bigwedge^3\g_i)^{\g_i}$ is one--dimensional
and generated by $[\Omega_{12}^i,\Omega_{23}^i]$, where
$\Omega^i\in(\g_i\otimes\g_i)^{\g_i}$ is the symmetric
element corresponding to the Killing form of $\g_i$, $\Phi$
is a non--degenerate solution of the pentagon equation iff
\begin{equation}\label{eq:non deg}
\Alt_3(\varphi)=\frac{1}{6}[\Omega_{12},\Omega_{23}]
\end{equation}
where $\Omega\in(\g\otimes\g)^{\g}$ corresponds to a
non--degenerate, ad--invariant, symmetric bilinear
form $(\cdot,\cdot)$ on $\g$.\\

Let now $D\subseteq\Dg$ be a subdiagram, and set
$$
\Phi_D=
\ol{\pi}^3_D(\Phi)\in
1+\hbar^2\left(\UgD^{\otimes 3}\fml\right)^{\gD}
$$
where $\ol{\pi}^3_D$ is the generalised \HC homomorphism
defined in \S \ref{ss:HC}. By proposition \ref{pr:HC face},
$\Phi_D$ satisfies the pentagon equation. Note that $\Phi_D$
is non--degenerate if $\Phi$ is.
 
\begin{theorem}\label{th:coh existence}
If $\Phi$ is non--degenerate, there exists an element
$$
F\in 1^{\otimes 2}+\hbar\left(\Ug^{\otimes 2}\fml\right)^{\lD}
$$
such that
\begin{equation}\label{eq:coh equation}
(\Phi)_F=\Phi_D
\quad\text{and}\quad
\ol{\pi}^2_D(F)=1\otimes 1
\end{equation}
Modulo $\hbar^2$, one has
$$
F=1^{\otimes 2}+\hbar(\rg-\rD)
$$
where $\rg,\rD$ are the standard solutions of the MCYBE for
$\g$ and $\gD$ corresponding to $(\cdot,\cdot)$. If $\Phi$
satisfies in addition
\begin{equation}\label{eq:constraint}
\Phi^{321}=\Phi^{-1}
\quad\text{and}\quad
\Phi^{\Theta}=\Phi
\end{equation}
where $\Theta\in\Aut(\g)$ is an involution acting as $-1$ on
$\h$, then $F$ may be chosen such that
\begin{equation}\label{eq:theta}
F^{\Theta}=F^{21}
\end{equation}
\end{theorem}

\remark If $\Phi$ is an associator, that is satisfies in addition $\id\otimes
\veps\otimes\id(\Phi)=1^{\otimes 2}$, where $\veps:\Ug\rightarrow\IC$ is
the counit, then so is $\Phi_D$ since one checks that $\id^{\otimes i}\otimes
\veps\otimes\id^{\otimes(k-i-1)}\circ\pi^k_D=\pi^{k-1}_D\circ\id^{\otimes i}
\otimes\veps\otimes\id^{\otimes (k-i-1)}$. In this case, it follows from 
\eqref{eq:coh equation} that $F$ is a twist, \ie satisfies $\veps\otimes
\id(F)=1=\id\otimes\veps(F)$.

\subsection{}

The proof of theorem \ref{th:coh existence} is given in
\S \ref{sss:ladder}--\S \ref{sss:secondary cocycle}. It closely
follows the argument of Donin--Shnider \cite[\S 3]{DS}
where theorem \ref{th:coh existence} is proved, under
the additional assumption \eqref{eq:constraint}, in the
case $D=\emptyset$. The reader familiar with Donin and
Shnider's argument will readily recognize that a relevant
difference is that the cohomology group
\begin{equation}\label{eq:old cohomology}
H^3(\bigwedge\g;\schl\rg,\cdot\schr)
\cong
\bigwedge^3\h
\end{equation}
which governs the secondary obstructions theory in \cite{DS}
is replaced by the group
\begin{equation}\label{eq:new cohomology}
H^3((\bigwedge\g)^{\gD};\schl\rg-\rD,\cdot\schr)
\cong
\bigwedge^3\cD\bigoplus(\bigwedge^3\gD)^{\gD}
\end{equation}
which was computed in \S \ref{ss:CE}. Another
significant difference is that the possibility of
laddering down from $\Dg$ to $D$ through intermediate
diagrams, as explained in \S \ref{sss:ladder}, allows
in effect to assume that $\cD$ is at most
two--dimensional, thus killing the first component
of the secondary obstruction in \eqref
{eq:new cohomology} and rendering the assumption
\eqref{eq:constraint} unnecessary to prove the
existence of $F$.

\subsection{}\label{sss:ladder}

Although we will only use this from \S \ref{sss:mod h^n} onwards,
note that we may assume that $|\Dg\setminus D|\leq 2$. Indeed,
assume theorem \ref{th:coh existence} proved in this case and let
$$\Dg=D_1\supset D_2
\supset\cdots\supset
D_{m-1}\supset D_{m}=D$$
be a nested chain of diagrams such that $|D_j\setminus D_{j+1}|
\leq 2$. For any pair $D''\subseteq D'\subseteq\Dg$, denote by
$\c_{D'',D'}\subset\h$ the span of the fundamental coweights
$\cow{k}$, with $k$ such that $\alpha_k\in D'\setminus D''$ and by
$$
\ol{\pi}_{D'',D'}^k:
(\Ug_{D'}\ok)^{\c_{D'',D'}}
\rightarrow
\Ug_{D''}\ok
$$
the corresponding generalised \HC homomorphism. Set $\Phi_1
=\Phi$ and, for $j=1\ldots m-1$
$$
\Phi_{j+1}=
\ol{\pi}^3_{D_{j+1},D_j}(\Phi_j)=
\ol{\pi}^3_{D_{j+1},\Dg}(\Phi)
$$
so that $\Phi_{m}=\Phi_D$. Let
$$
F_j\in 1^{\otimes 2}+\hbar(\Ug_{D_j}^{\otimes 2}\fml)^{\g_{D_{j+1}}}
$$
be such that
$$
(\Phi_j)_{F_j}=\Phi_{j+1}
\quad\text{and}\quad
\ol{\pi}^2_{D_{j+1},D_j}(F_j)=1^{\otimes 2}
$$
then
$$
F=
F_{m-1}\cdots F_1
$$
is readily seen to satisfy \eqref{eq:coh equation}.
Note that if $\Phi$ satisfies in addition \eqref
{eq:constraint} then so does each $\Phi_j$ since
$\ol{\pi}_D^k$ is equivariant for the action of
the symmetric group $\SS_k$ and of $\Theta$. In
that case, choosing each $F_j$ such that $F_j
^{\Theta}=F_j^{21}$ yields an $F$ which
satisfies $F^{\Theta}=F^{21}$.

\subsection{}\label{sss:mod h^2}

We begin by solving equation \eqref{eq:coh equation} mod $\hbar^2$.
Let $f\in(\Ug^{\otimes 2})^{\lD}$ and set $F=1^{\otimes 2}+\hbar f$.
Since $\Phi$ and $\Phi_D$ are equal to $1^{\otimes 3}$ mod $\hbar
^2$, the coefficient of $\hbar$ in $(\Phi)_F-\Phi_D$ is
$$1\otimes f+\id\otimes\Delta(f)-\Delta\otimes\id(f)-f\otimes 1
=
\dH f$$
Thus, $F$ is a solution of \eqref{eq:coh equation} mod $\hbar^2$ if,
and only if, $f$ is a \Ho 2--cocycle such that $\ol{\pi}_D^2f=0$.

\subsection{}\label{sss:general step}

Let now $n\geq 1$ and let
$$F=1^{\otimes 2}+\hbar f+\cdots+\hbar^{n}f_{n}
\in
1^{\otimes 2}+\hbar(\Ug^{\otimes 2}\fml)^{\lD}$$
be a solution of \eqref{eq:coh equation} mod $\hbar^{n+1}$.
We shall derive below a necessary and sufficient condition
for \eqref{eq:coh equation} to possess a solution mod $\hbar
^{n+2}$ of the form $\wt{F}=F+\hbar^{n+1}f_{n+1}$ where
$$f_{n+1}\in(\Ug^{\otimes 2})^{\lD}
\qquad\text{satisfies}\qquad
\ol{\pi}_D^2(f_{n+1})=0$$
Define $\xi\in\Ug^{\otimes 3}$ by
\begin{equation}\label{eq:xi}
1\otimes F\cdot\id\otimes\Delta(F)\cdot\Phi-
\Phi_D\cdot F\otimes 1\cdot\Delta\otimes\id(F)=
\hbar^{n+1}\xi
\thickspace\mod\hbar^{n+2}
\end{equation}
Then, $\wt F$ is a solution of \eqref{eq:coh equation} mod
$\hbar^{n+2}$ if, and only if $\dH f_{n+1}=-\xi$.

\begin{lemma}\label{le:deviation}
The element $\xi$ is invariant under $\lD$ and satisfies
$$
\dH\xi=0
\quad\text{and}\quad
\ol{\pi}_D^3(\xi)=0
$$
\end{lemma}
\proof The invariance of $\xi$ under $\lD$ follows from
that of $\Phi,\Phi_D$ and $F$. Since $F=1^{\otimes 2}$
mod $\hbar$,
$$\hbar^{n+1}\xi=(\Phi)_F-\Phi_D\mod\hbar^{n+2}$$
Since $F$ is invariant under $\gD$, the restriction of
$$\Delta_F(\cdot)=F\cdot\Delta(\cdot)\cdot F^{-1}$$
to $\UgD$ is equal to $\Delta$ and $\Phi_D$ satisfies
the pentagon equation with respect to $\Delta_F$. Since
this is also the case of $(\Phi)_F$, we find, working
mod $\hbar^{n+2}$, that
$$
0=
\Pent_{\Delta_F}((\Phi)_F)=
\Pent_{\Delta_F}(\Phi_D)+\hbar^{n+1}\dH\xi=
\hbar^{n+1}\dH\xi
$$
where, for any $\Psi\in\Ug^{\otimes 3}$ and map $\wt{\Delta}
:\Ug\rightarrow\Ug^{\otimes 2}$,
$$
\Pent_{\wt{\Delta}}(\Psi)=
1\otimes\Psi\cdot
\id\otimes\wt{\Delta}\otimes\id(\Psi)\cdot
\Psi\otimes 1
-
\id^{\otimes 2}\otimes\wt{\Delta}(\Psi)
\cdot
\wt{\Delta}\otimes\id^{\otimes 2}(\Psi)
$$
Finally, from $\ol{\pi}_D^2(F)=1^{\otimes 2}$ and
$\Phi_D=\ol{\pi}_D^3(\Phi)$, we get, using proposition 
\ref{pr:HC face} that
$$
\hbar^{n+1}\ol{\pi}_D^3\xi=
(\ol{\pi}_D^3(\Phi))_{\ol{\pi}_D^2(F)}-\Phi_D=
0
$$
\halmos

\begin{lemma}\label{le:theta deviation}
If $\Phi$ and $F$ satisfy \eqref{eq:constraint} and \eqref
{eq:theta} respectively, then
\begin{equation}\label{eq:theta deviation}
\xi^{\Theta}=-\xi^{321}
\end{equation}
\end{lemma}
\proof We have, working mod $\hbar^{n+2}$,
\begin{equation*}
\begin{split}
\Phi_D^{\Theta}+\hbar^{n+1}\xi^{\Theta}
&=
1\otimes F^{21}\cdot\id\otimes\Delta(F^{21})
\cdot\Phi\cdot
\Delta\otimes\id((F^{21})^{-1})\otimes(F^{21})^{-1}\otimes 1 \\
&=
\left(
F\otimes 1\cdot\Delta\otimes\id(F)
\cdot\Phi^{321}\cdot
\id\otimes\Delta(F^{-1})\cdot 1\otimes F^{-1}
\right)^{321}\\
&=
\left( (
1\otimes F\cdot\id\otimes\Delta(F)
\cdot\Phi\cdot
\Delta\otimes\id(F^{-1})\cdot F^{-1}\otimes 1)^{-1}\right)^{321}\\
&=
\left((\Phi_D+\hbar^{n+1}\xi)^{-1}\right)^{321}\\
&=
(\Phi_D^{-1}-\hbar^{n+1}\xi)^{321}
\end{split}
\end{equation*}
whence \eqref{eq:theta deviation} since $\Phi_D=\ol{\pi}_
{D}^3(\Phi)$ satisfies \eqref{eq:constraint} \halmos

\begin{corollary}\label{co:primary extension}
The element $F$ may be extended to a solution of \eqref
{eq:coh equation} mod $\hbar^{n+2}$ if, and only if $
\Alt_3\xi=0$. If in addition $\Phi,F$ satisfy \eqref
{eq:constraint} and \eqref{eq:theta} respectively, the
extension may be chosen so as to satisfy \eqref{eq:theta}.
\end{corollary}
\proof 
$\Alt_3\xi=0$ if, and only if $\xi=\dH g$ for some $g\in\Ug^
{\otimes 2}$, which may then be chosen invariant under
$\lD$. By corollary \ref{co:HC Ho}, we have
$$
0=
\ol{\pi}_D^3\xi =
\ol{\pi}_D^3\dH g=
\dH\ol{\pi}_D^2 g
$$
so that, setting $f_{n+1}=-(g-\ol{\pi}_D^2(g))$ we
have
$$
\dH f_{n+1}=-\xi
\qquad\text{and}\qquad
\ol{\pi}_D^2f_{n+1}=0
$$
and $F+\hbar^{n+1}f_{n+1}$ is a solution of \eqref
{eq:coh equation} mod $\hbar^{n+2}$. If $\Phi,F$ satisfy
\eqref{eq:constraint} and \eqref{eq:theta} respectively,
then, by lemma \ref{le:theta deviation}
$$
 \dH f_{n+1}^{\Theta}=
-\xi^{\Theta}=
 \xi^{321}=
-(\dH f_{n+1})^{321}=
 \dH f_{n+1}^{21}
$$
so that $f_{n+1}'=1/2(f_{n+1}+(f_{n+1}^{21})^{\Theta})$
satisfies
$$
\dH f_{n+1}'=\xi,
\quad
\ol{\pi}_D^2f_{n+1}'=0
\quad\text{and}\quad
(f_{n+1}')^{\Theta}=
(f_{n+1}')^{21}
$$
and $F+\hbar^{n+1}f_{n+1}'$ solves \eqref{eq:coh equation}
mod $\hbar^{n+2}$ and satisfies \eqref{eq:theta} \halmos
 
\subsection{}\label{sss:mod h^3}

We consider first the case $n=1$ so that $F=1+\hbar f$ where
$f\in(\Ug^{\otimes 2})^{\lD}$ is a \Ho 2--cocycle such that
$\ol{\pi}_D^2(f)=0$. By lemma \ref{le:gauge ext} below,
adding a 2--coboundary to $f$ does not affect the extendability
of $F$ to a solution mod $\hbar^3$. We may therefore assume
that $f\in(\bigwedge^2\g)^{\lD}$. In this case, since $\Phi$
and $\Phi_D$ are equal to $1^{\otimes 3}$ mod $\hbar^2$,
we get
$$
\xi=
\varphi-\varphi_D+
f^{23}(f^{12}+f^{13})-f^{12}(f^{13}+f^{23})
$$
where $\Phi=1+\hbar^2\varphi$ mod $\hbar^3$ and $\varphi
_D=\ol{\pi}_D^3\varphi$. Thus, $F$ extends to a solution
mod $\hbar^3$ if, and only if,
$$
\Alt_3(\varphi)-\Alt_3(\varphi_D)
=
\Alt_3(f^{12}(f^{13}+f^{23})-f^{23}(f^{12}+f^{13}))
$$
We shall need the following

\begin{lemma}\label{le:alt=sch}
For any $f,\chi\in\bigwedge^2\g$, one has
\begin{equation}\label{eq:alt=sch}
\begin{split}
&
\Alt_3\left(
 f^{12}(\chi^{13}+\chi^{23})+\chi^{12}(f^{13}+f^{23})
-f^{23}(\chi^{12}+\chi^{13})-\chi^{23}(f^{12}+f^{13})
\right)\\
=&
\schl f,\chi\schr
\end{split}
\end{equation}

where $\schl\cdot,\cdot\schr$ is the Schouten bracket
\eqref{eq:schouten}.
\end{lemma}
\proof Since
\begin{gather*}
(f^{12}(\chi^{13}+\chi^{23}))^{(1\medspace 3)}
=
f^{23}(\chi^{13}+\chi^{12}),\\
(f^{12}\chi^{23})^{(1\medspace 2)}=-f^{12}\chi^{13}
\aand
(\chi^{12}f^{13})^{(2\medspace 3)}=\chi^{13}f^{12}
\end{gather*}
the left--hand side of \eqref{eq:alt=sch} is equal to
\begin{equation*}
\begin{split}
2\Alt_3\left(
 f^{12}(\chi^{13}+\chi^{23})+\chi^{12}(f^{13}+f^{23})
 \right)
&=
4\Alt_3(f^{12}\chi^{13}+\chi^{12}f^{13})\\
&=
4\Alt_3([f^{12},\chi^{13}])\\
&=
\frac{2}{3}\YB(f,\chi)\\
&=
\schl f,\chi\schr
\end{split}
\end{equation*}
where we used \eqref{eq:Alt=YB} and lemma \ref{le:YB=sch}
\halmos

\begin{corollary}\label{co:mod h^3}
Let $f\in(\bigwedge^2\g)^{\lD}$ be such that $\ol{\pi}_D^2(f)=0$.
Then, the element
$$F=1+\hbar f$$
extends to a solution $\wt F$ of \eqref{eq:coh equation} mod
$\hbar^3$ if, and only if,
$$\half{1}\schl f,f\schr=
\Alt_3\varphi-\Alt_3\varphi_D$$
In that case, and provided \eqref{eq:constraint} holds, $\wt F$
may be chosen so as to satisfy \eqref{eq:theta}.
\end{corollary}

Let $(\cdot,\cdot)$ be the non--degenerate, ad--invariant, symmetric
bilinear form on $\g$ such that
\begin{equation}\label{eq:phi Omega}
\Alt_3(\varphi)=
\frac{1}{6}[\Omega_{12},\Omega_{23}]
\end{equation}
and $\rg,\rD$ be the standard solutions of the MCYBE determined
by $(\cdot,\cdot)$ and its restriction to $\gD$ respectively, so that
$\ol{\pi}_D^2(\rg)=\rD$. We henceforth set
\begin{equation}\label{eq:f}
f=\rg-\rD\in(\bigwedge^2\g)^{\lD}
\end{equation}
By corollary \ref{co:mod h^3}, $F=1+\hbar f$ extends to a solution
of \eqref{eq:coh equation} mod $\hbar^3$ which, in addition, satisfies
\eqref{eq:theta} if \eqref{eq:constraint} holds. Indeed, by proposition
\ref{pr:differential}, lemma \ref{le:YB=sch} and proposition
\ref{pr:MCYBE}, we have
\begin{multline*}
\schl f,f\schr
=
\schl\rg,\rg\schr-\schl\rD,\rD\schr\\
=
\frac{2}{3}\left(\YB(\rg,\rg)-\YB(\rD,\rD)\right)=
\frac{1}{3}\left(
[\Omega_{12},\Omega_{23}]-[\Omega_{12}^D,\Omega_{23}^D]\right)
\end{multline*}

\subsection{}\label{sss:mod h^n}

Assume now $n\geq 2$ and let
$$F=1+\hbar f+\hbar^2f_2+\cdots+\hbar^{n}f_n$$
be a solution of \eqref{eq:coh equation} mod $\hbar^{n+1}$.
Let $\xi=\xi(f;f_2,\ldots,f_{n})\in\Ug^{\otimes 3}$ be
given by \eqref{eq:xi}. By \S \ref{sss:general step}, $\xi$
is a \Ho 3--cocycle and $F$ extends to a solution mod $\hbar
^{n+2}$ if, and only if $\xi$ is a coboundary. This, however
need not be the case. We note none--the--less that if $\chi
\in(\Ug^{\otimes 2})^{\lD}$ satisfies
$$
\dH\chi=0
\quad\text{and}\quad
\ol{\pi}_D^2(\chi)=0
$$
then $F+\hbar^{n}\chi$ is also a solution of \eqref
{eq:coh equation} mod $\hbar^{n+1}$ which could admit
an extension mod $\hbar^{n+2}$. By the following result,
the extendability of $F+\hbar^{n}\chi$ only depends
upon the \Ho cohomology class of $\chi$.

\begin{lemma}\label{le:gauge ext}
If $\chi$ is a \Ho 2--coboundary, then the element $F+\hbar
^{n}\chi$ can be extended to a solution mod $\hbar^{n+2}$
of \eqref{eq:coh equation} if, and only if $F$ can.
\end{lemma}
\proof It suffices to prove one implication since $F=
(F+\hbar^{n}\chi)-\hbar^{n}\chi$. Write $\chi=\dH g$
with $g\in\Ug^{\lD}$. By corollary \ref{co:HC Ho},
$$
0=\ol{\pi}_D^2(\chi)=\dH\ol{\pi}_D^1g
$$
so that $\chi=\dH g'$ where $g'=(1-\ol{\pi}_D^1)g$
is invariant under $\lD$ and lies in the kernel of $\ol
{\pi}_D^1$. Let $\wt{F}=F+\hbar^{n+1}f_{n+1}$ be a
solution mod $\hbar^{n+2}$ of \eqref{eq:coh equation}.
Then
$$
\wt{F}'=
u^{\otimes 2}\cdot\wt{F}\cdot\Delta(u)^{-1}
$$
where $u=1+\hbar^{n}g¹$, is equal to $F+\hbar^{n}\chi$
mod $\hbar^{n+1}$ and solves \eqref{eq:coh equation}
mod $\hbar^{n+2}$ since 
\begin{equation*}
\begin{split}
(\Phi)_{\wt{F}'}
&=
u^{\otimes 3}\cdot
1\otimes\wt{F}\cdot\id\otimes\Delta(\wt{F})\cdot
\id\otimes\Delta(\Delta(u))^{-1}\cdot
\Phi\\
&\phantom{=}
\cdot
\Delta\otimes\id(\Delta(u))\cdot
\Delta\otimes\id(\wt{F})\cdot\wt{F}\otimes 1\cdot
(u^{\otimes 3})^{-1}\\
&=
u^{\otimes 3}\cdot
\Phi_D\cdot
(u^{\otimes 3})^{-1}\\
&=
\Phi_D
\end{split}
\end{equation*}
where the first equality follows from the $\g$--invariance
of $\Phi$ and the last from the $\gD$--invariance of $u$
\halmos\\

We may therefore assume that $\chi$ lies in $(\bigwedge^2
\g)^{\lD}$. We then note that, for $n\geq 2$,
\begin{equation*}
\begin{split}
\xi(f;f_2,\ldots,f_{n}+\chi)
&=
\xi(f;f_2,\ldots,f_{n})
+
f^{23}(\chi^{12}+\chi^{13})+\chi^{23}(f^{12}+f^{13})\\
&\phantom{=}
-f^{12}(\chi^{13}+\chi^{23})-\chi^{12}(f^{13}+f^{23})
\end{split}
\end{equation*}
so that $F+\hbar^{n}\chi$ possesses an extension mod $\hbar
^{n+2}$ if, and only if,
$$\Alt_3(\xi(f;f_2,\ldots,f_{n}))=
\schl f,\chi\schr$$
where we used lemma \ref{le:alt=sch}

\begin{proposition}\label{pr:sec cocycle}
The element $\wt{\xi}=\Alt_3(\xi)\in(\bigwedge^3\g)^{\lD}$
satisfies
$$
\schl f,\wt{\xi}\schr=0
\quad\text{and}\quad
\ol{\pi}_D^2(\wt{\xi})=0
$$
\end{proposition}

We defer the proof of proposition \ref{pr:sec cocycle}
to \S \ref{sss:secondary cocycle} in order to conclude
the proof of theorem \ref{th:coh existence}. By proposition
\ref{pr:sec cocycle}, $\wt\xi$ is a 3--cocycle in $\left
((\extg{})^{\lD},\schl f,\cdot\schr\right)$ and we must
show that it is a 3--coboundary. By theorem \ref{th:sch coh}
\begin{equation}\label{eq:almost}
\wt\xi=\schl f,\chi\schr+\eta
\end{equation}
for some $\chi\in(\extg{2})^{\lD}$ and
$$
\eta\in
(\bigwedge^3\lD)^{\lD}
=
\bigwedge^3\cD\oplus(\bigwedge^3\gD)^{\gD}
=
(\bigwedge^3\gD)^{\gD}
$$
where the first equality follows from the fact that 
$(\bigwedge^i\gD)^{\gD}=0$ for $i=1,2$ and the second
from the assumption that $|\Dg\setminus D|\leq 2$ so that
$\cD$ is at most two--dimensional. Applying $\ol{\pi}_D^3$
to both sides of \eqref{eq:almost}, we find, since $\ol{\pi}_
D^2(f)=0$, that
$$0=
\ol{\pi}_D^3(\schl f,\chi\schr+\eta)=
\eta$$

Note that if $\Phi,F$ satisfy \eqref{eq:constraint} and \eqref
{eq:theta} respectively, then by lemma \ref{le:theta deviation}
\begin{equation}\label{eq:theta xi}
\wt\xi^{\Theta}=
\Alt_3(\xi^{\Theta})=
\Alt_3(-\xi^{321})=
\wt\xi
\end{equation}
Since $f^{\Theta}=-f$, this implies
$$
\schl f,\chi\schr=
\schl f,\chi\schr^{\Theta}=
-\schl f,\chi^{\Theta}\schr
$$
so that $\chi'=1/2(\chi-\chi^{\Theta})$ satisfies $\chi'^{\Theta}=
-\chi'=\chi'^{21}$ and $\schl f,\chi'\schr=\wt\xi$ and $F+\hbar^{n}
\chi'$ is a solution of \eqref{eq:coh equation} mod $\hbar^{n+1}$
possessing an extension to a solution mod $\hbar^{n+2}$ satisfying
\eqref{eq:theta}. This completes the proof of theorem \ref
{th:coh existence} \halmos\\

\remark If $\Phi$ satisfies \eqref{eq:constraint}, it is not
necessary to ladder down, \ie assume that $|\Dg\setminus D|
\leq 2$. Indeed, by theorem \ref{th:sch coh}, there exist
unique elements $u\in\bigwedge^3\cD$, $v\in(\bigwedge^3
\gD)^{\gD}$ and a $\chi\in(\bigwedge^2\g)^{\lD}$ such that
$$
\wt\xi
=u+v+\schl f,\chi\schr
$$
Since $\wt\xi$ and $u$ are killed by $\ol{\pi}_D^3$, $v=0$.
Applying $\Theta$ to both sides and using \eqref{eq:theta xi},
we find that $u^{\Theta}=u$. This however implies that $u=0$
since $\Theta$ acts by $-1$ on $\bigwedge^3\h\supseteq\bigwedge
^3\cD$.

\subsection{Proof of proposition \ref{pr:sec cocycle}}
\label{sss:secondary cocycle}

We begin with some preliminary lemmas. Let $\wt{\Delta}:\Ug
\rightarrow\Ug^{\otimes 2}$ be a linear map. Let 
$$
\xi=\xi_1\otimes\cdots\otimes\xi_k\in\Ug^{\otimes k}
$$
$i\leq k$, and write
$$
\wt{\Delta}(\xi_i)=\sum_{a}\xi_{i,a}'\otimes\xi_{i,a}''
$$
For any enumeration $j_1,\ldots,j_{k+1}$ of $[1,k+1]$, we
set
$$
\xi^{j_1,\ldots,j_{i-1},j_ij_{i+1},j_{i+2},\ldots,j_{k+1}}
=
\sum_{a}\eta_{a}
$$
where $\eta_{a}\in\Ug^{\otimes(k+1)}$ is the decomposable
tensor with component $\xi_{\ell}$ in position $j_{\ell}$
if $\ell\leq i-1$ and in position $j_{\ell+1}$ if $\ell\geq
i+1$ and  components $\xi_{i,a}',\xi_{i,a}''$ in positions
$j_i$ and $j_{i+1}$ respectively. In other words,
\begin{equation*}
\begin{split}
\xi^{j_1,\ldots,j_{i-1},j_ij_{i+1},j_{i+2},\ldots,j_{k+1}}
&=
\sigma
\xi^{1,\ldots,i-1,i\thinspace i+1,i+2,\ldots,k+1}\\
&=
\sigma\circ
\id^{\otimes(i-1)}\otimes\wt{\Delta}\otimes\id^{\otimes(k-i)}
\medspace\xi
\end{split}
\end{equation*}
where $\sigma\in\SS_{k+1}$ is the permutation mapping $\ell$
to $i_{\ell}$.

\begin{lemma}\label{le:alt alt}
For any $\xi\in\extg{k}$, one has
\begin{multline}\label{eq:alt alt}
(k+1)\Alt_{k+1}\left(
\sum_{i=1}^k(-1)^i
\id^{\otimes(i-1)}\otimes\wt{\Delta}\otimes\id^{(k-i)}\xi
\right)\\
=
\sum_{1\leq a<b\leq k+1}(-1)^{a+b}\left(
(\Alt_k\xi)^{ab,1,\ldots,\wh{a},\ldots,\wh{b},\ldots,k+1}-
(\Alt_k\xi)^{ba,1,\ldots,\wh{a},\ldots,\wh{b},\ldots,k+1}
\right)
\end{multline}
\end{lemma}
\proof For any $i\in[1,k]$ and $a\neq b\in[1,k+1]$, set
$$
\SS_{k+1}^{i:a,b}=
\{\sigma\in\SS_{k+1}|
  \medspace\sigma(i)=a,\medspace\sigma(i+1)=b\}
$$
Then, the left--hand side of \eqref{eq:alt alt} is equal to
\begin{equation*}
\begin{split}
 &
\frac{1}{k!}
\sum_{\substack
{1\leq a<b\leq k+1\\[.6ex]
 1\leq i\leq k\\[.6ex]
 \sigma\in\SS_{k+1}^{i:a,b}}}
(-1)^i(-1)^{\sigma}
\left(
\xi^{\sigma(1),\ldots,\sigma(i-1),ab,
     \sigma(i+2),\ldots,\sigma(k+1)}
\right.\\[-4em]
&\phantom{======================}
-
\left.
\xi^{\sigma(1),\ldots,\sigma(i-1),ba,
     \sigma(i+2),\ldots,\sigma(k+1)} 
\right)\\[1.6em]
=&
\frac{1}{k!}
\sum_{\substack
{1\leq a<b\leq k+1\\[.6ex]
 1\leq i\leq k\\[.6ex]
 \sigma\in\SS_{k+1}^{i:a,b}}}
(-1)^i(-1)^{\sigma}
\left(
(\ol{\sigma}\xi)^{ab,1,\ldots,\wh{a},\ldots,\wh{b}\ldots,k+1}-
(\ol{\sigma}\xi)^{ba,1,\ldots,\wh{a},\ldots,\wh{b}\ldots,k+1}
\right)
\end{split}
\end{equation*}
where, for any $\sigma\in\SS_{k+1}^{i:a,b}$,
$$\ol{\sigma}\in\SS_k^i=
\{\tau\in\SS_k|\medspace\tau(i)=1\}$$
is the permutation determined by the commutativity of the 
following
diagram
$$
\begin{diagram}
[1,k]\setminus\{i\}&\rTo&[1,k+1]\setminus\{i,i+1\}\\
\dTo^{\ol{\sigma}} &	&\dTo_{\sigma}		  \\
[1,k]\setminus\{1\}&\rTo&[1,k+1]\setminus\{a,b\}
\end{diagram}
$$
where the horizontal arrow are the obvious monotone 
identifications. Noting that $\sigma\rightarrow\ol
{\sigma}$ is an isomorphism of $\SS_{k+1}^{i:a,b}$
onto $\SS_k^i$ and deferring for the time being
the proof that
\begin{equation}\label{eq:sign sign}
(-1)^i(-1)^{\sigma}=(-1)^{a+b}(-1)^{\ol{\sigma}}
\end{equation}
we see that the above is equal to 
$$
\frac{1}{k!}
\sum_{\substack
{1\leq a<b\leq k+1\\[.6ex]
 1\leq i\leq k\\[.6ex]
 \ol{\sigma}\in\SS_k^i}}
(-1)^{a+b}(-1)^{\ol{\sigma}}
\left(
(\ol{\sigma}\xi)^{ab,1,\ldots,\wh{a},\ldots,\wh{b}\ldots,k+1}-
(\ol{\sigma}\xi)^{ba,1,\ldots,\wh{a},\ldots,\wh{b}\ldots,k+1}
\right)
$$
and therefore to the right--hand side of \eqref{eq:alt alt}.
We turn now to the proof of \eqref{eq:sign sign}. Let $\ol{\ol
{\sigma}}\in\SS_{k-1}$ be the permutation determined by the
commutativity of
$$
\begin{diagram}
[1,k]\setminus\{i\}&\rTo&[1,k-1]                &\rTo&[1,k+1]\setminus\{i,i+1\}\\
\dTo^{\ol{\sigma}} &	&\dTo^{\ol{\ol{\sigma}}}&    &\dTo^{\sigma}	       \\
[1,k]\setminus\{1\}&\rTo&[1,k-1]                &\rTo&[1,k+1]\setminus\{a,b\}
\end{diagram}
$$
where the horizontal arrows are the obvious monotone identifications.
Since $(-1)^{\ol{\sigma}}=(-1)^{\ol{\ol{\sigma}}}\cdot(-1)^{i-1}$,
it suffices to prove that $(-1)^{\sigma}=(-1)^{\ol{\ol{\sigma}}}
(-1)^{a+b-1}$. This clearly holds if $a=1$ and $b=2$. In the
general case, letting $\tau\in\SS_{k+1}$ be the unique permutation
such that $\tau$ is increasing on $[1,k+1]\setminus\{a,b\}$, $\tau
(a)=1$ and $\tau(b)=2$, so that $(-1)^{\tau}=(-1)^{a+b-1}$, and
noting that $\ol{\ol{\tau\circ\sigma}}=\ol{\ol{\sigma}}$, we see
that
$$
(-1)^{\sigma}=
(-1)^{a+b-1}(-1)^{\tau\circ\sigma}=
(-1)^{a+b-1}(-1)^{\ol{\ol{\tau\circ\sigma}}}=
(-1)^{a+b-1}(-1)^{\ol{\ol{\sigma}}}
$$
\halmos

\begin{lemma}\label{le:1 k}
For any $Y,X_1,\ldots,X_k\in\g$, one has
\begin{multline}\label{eq:1 k}
Y\wedge X_1\wedge\ldots\wedge X_k\\
=
\frac{1}{(k+1)!}
\sum_{\substack{1\leq i\leq k+1\\[.6ex]\tau\in\SS_k}}
(-1)^{i-1}(-1)^{\tau}
X_{\tau(1)}\otimes\cdots\otimes X_{\tau(i-1)}\otimes
Y\otimes
X_{\tau(i)}\otimes\cdots\otimes X_{\tau(k)}
\end{multline}
\end{lemma}
\proof Set $Z_1=Y$ and $Z_j=X_{j-1}$ for $j=2\ldots k+1$.
By definition, $(k+1)!$ times the \lhs of \eqref{eq:1 k} is
equal to
\begin{equation*}
\begin{split}
 &
\sum_{\tau\in\SS_{k+1}}(-1)^{\tau}
Z_{\tau(1)}\otimes\cdots\otimes Z_{\tau(k+1)}\\
=&
\sum_{\substack{1\leq j\leq k+1\\[.6ex]\tau\in\SS_{k+1}:\tau(j)=1}}
(-1)^{\tau}
X_{\tau(1)-1}\otimes\cdots\otimes X_{\tau(j-1)-1}
\otimes Y\otimes X_{\tau(j+1)-1}\otimes\cdots\\[-2em]
&\phantom{================================}
\cdots\otimes X_{\tau(k+1)-1}
\end{split}
\end{equation*}
For any $\tau\in\SS_{k+1}$ such that $\tau(j)=1$, let $\ol{\tau}
\in\SS_k$ be the permutation determined by the commutativity
of 
$$
\begin{diagram}
[1,k]		&\rTo&[1,k+1]\setminus\{j\}\\
\dTo^{\ol{\tau}}&    &\dTo_{\tau}          \\
[1,k]		&\rTo&[1,k+1]\setminus\{1\}
\end{diagram}
$$
Then, $(-1)^{\tau}=(-1)^{\ol{\tau}}(-1)^{j-1}$ and the above is
equal to 
$$
\sum_{\substack{1\leq j\leq k+1\\[.6ex]\tau\in\SS_{k+1}:\tau(j)=1}}
(-1)^{j-1}(-1)^{\ol{\tau}}
X_{\ol{\tau}(1)}\otimes\cdots\otimes X_{\ol{\tau}(j-1)}\otimes
Y\otimes 
X_{\ol{\tau}(j)}\otimes\cdots\otimes X_{\ol{\tau}(k)}
$$
which proves \eqref{eq:1 k} \halmos

\begin{lemma}\label{le:alt sch}
For any $f\in\extg{2}$ and $\eta\in\extg{k}$, one has
\begin{equation}\label{eq:alt sch}
\sum_{1\leq a<b\leq k+1}(-1)^{a+b}
[f^{ab},
\eta^{a,1,\ldots,\wh{a},\ldots,\wh{b},\ldots,k+1}+
\eta^{b,1,\ldots,\wh{a},\ldots,\wh{b},\ldots,k+1}]
=
-\frac{k+1}{2}\schl f,\eta\schr
\end{equation}
\end{lemma}
\proof We may assume that $f,\eta$ are of the form
\begin{gather*}
f=
f_1\wedge f_2=
\frac{1}{2}(f_1\otimes f_2-f_2\otimes f_1)\\
\eta=
\eta_1\wedge\ldots\wedge\eta_k=
\frac{1}{k!}\sum_{\sigma\in\SS_k}(-1)^{\sigma}
\eta_{\sigma(1)}\otimes\cdots\otimes\eta_{\sigma(k)}
\end{gather*}
The left--hand side of \eqref{eq:alt sch} is then equal to
\begin{equation}\label{eq:summands}
\begin{split}
 &
\frac{1}{2k!}
\sum_{\substack{1\leq a<b\leq k+1\\ \sigma\in\SS_k}}
(-1)^{a+b}(-1)^{\sigma}\\
&\phantom{=-}
[f_1^{a}f_2^{b}-f_2^{a}f_1^{b},
(\eta_{\sigma(1)}\otimes\cdots\otimes\eta_{\sigma(k)})^
{a,1,\ldots,\wh{a},\ldots,\wh{b},\ldots,k+1}\\
&\phantom{==================}
+
(\eta_{\sigma(1)}\otimes\cdots\otimes\eta_{\sigma(k)})^
{b,1,\ldots,\wh{a},\ldots,\wh{b},\ldots,k+1}]\\
=&
\frac{1}{2k!}
\sum_{\substack{1\leq a<b\leq k+1\\ \sigma\in\SS_k}}
(-1)^{a+b}(-1)^{\sigma}\\
&\phantom{=-}
\eta_{\sigma(2)}\otimes\cdots\otimes\eta_{\sigma(a)}\otimes
[f_1,\eta_{\sigma(1)}]\otimes
\eta_{\sigma(a+1)}\otimes\cdots\\
&\phantom{===============}\cdots
\otimes\eta_{\sigma(b-1)}\otimes
f_2\otimes
\eta_{\sigma(b)}\otimes\cdots\otimes\eta_{\sigma(k)}\\[.5em]
&\phantom{=}
-
\eta_{\sigma(2)}\otimes\cdots\otimes\eta_{\sigma(a)}\otimes
f_2\otimes
\eta_{\sigma(a+1)}\otimes\cdots\\
&\phantom{===============}\cdots
\otimes\eta_{\sigma(b-1)}\otimes
[f_1,\eta_{\sigma(1)}]\otimes
\eta_{\sigma(b)}\otimes\cdots\otimes\eta_{\sigma(k)}\\[.5em]
&\phantom{=}
-
\eta_{\sigma(2)}\otimes\cdots\otimes\eta_{\sigma(a)}\otimes
[f_2,\eta_{\sigma(1)}]\otimes
\eta_{\sigma(a+1)}\otimes\cdots\\
&\phantom{===============}\cdots
\otimes\eta_{\sigma(b-1)}\otimes
f_1\otimes
\eta_{\sigma(b)}\otimes\cdots\otimes\eta_{\sigma(k)}\\[.5em]
&\phantom{=}
+
\eta_{\sigma(2)}\otimes\cdots\otimes\eta_{\sigma(a)}\otimes
f_1\otimes
\eta_{\sigma(a+1)}\otimes\cdots\\
&\phantom{===============}\cdots
\otimes\eta_{\sigma(b-1)}\otimes
[f_2,\eta_{\sigma(1)}]\otimes
\eta_{\sigma(b)}\otimes\cdots\otimes\eta_{\sigma(k)}
\end{split}
\end{equation}
Setting $\sigma'=\sigma\circ(1\cdots a)$ in the first summand
and $\sigma'=\sigma\circ(1\cdots b-1)$ in the second, we see
that their sum is equal to
\begin{equation*}
\begin{split}
 &
\frac{1}{2k!}
\sum_{\substack{1\leq a<b\leq k+1\\ \sigma\in\SS_k}}
(-1)^{b-1}(-1)^{\sigma}\\
&\phantom{=-}
\eta_{\sigma(1)}\otimes\cdots\otimes\eta_{\sigma(a-1)}\otimes
[f_1,\eta_{\sigma(a)}]
\otimes\eta_{\sigma(a+1)}\otimes\cdots\\
&\phantom{===================}\cdots
\otimes\eta_{\sigma(b-1)}\otimes f_2\otimes
\eta_{\sigma(b)}\otimes\cdots\otimes\eta_{\sigma(k)}\\
-&
\frac{1}{2k!}
\sum_{\substack{1\leq a<b\leq k+1\\ \sigma\in\SS_k}}
(-1)^{a}(-1)^{\sigma}\\
&\phantom{=-}
\eta_{\sigma(1)}\otimes\cdots\otimes\eta_{\sigma(a-1)}
\otimes f_2\otimes
\eta_{\sigma(a)}\otimes\cdots\\
&\phantom{==============}\cdots
\otimes\eta_{\sigma(b-2)}\otimes
[f_1,\eta_{\sigma(b-1)}]\otimes
\eta_{\sigma(b)}\otimes\cdots\otimes\eta_{\sigma(k)}
\end{split}
\end{equation*}
and therefore to
\begin{equation*}
\begin{split}
&
\frac{1}{2k!}
\sum_{\substack{1\leq a\neq b\leq k\\\sigma\in\SS_k}}
(-1)^{b-1}(-1)^{\sigma}\\
&\phantom{=-}
\eta_{\sigma(1)}\otimes\cdots\otimes\eta_{\sigma(a-1)}\otimes
[f_1,\eta_{\sigma(a)}]\otimes
\eta_{\sigma(a+1)}\otimes\cdots\\
&\phantom{===================}
\cdots\otimes\eta_{\sigma(b-1)}\otimes
f_2\otimes
\eta_{\sigma(b)}\otimes\cdots\otimes\eta_{\sigma(k)}\\[.5em]
=&
\frac{k+1}{2}f_2\wedge
(\ad(f_1)\medspace\eta)
\end{split}
\end{equation*}
where we used lemma \ref{le:1 k}. Similarly, the sum of the
last two summands in \eqref{eq:summands} is equal to
$$
-\frac{k+1}{2}f_1\wedge
(\ad(f_2)\medspace\eta)
$$
Thus, the left--hand side of \eqref{eq:alt sch} is equal to
$$
\frac{k+1}{2}\left(
f_2\wedge\ad(f_1)\eta-
f_1\wedge\ad(f_2)\eta\right)
=
-\frac{k+1}{2}\schl f_1\wedge f_2,\eta\schr
$$ 
as claimed \halmos\\

{\sc Proof of proposition \ref{pr:sec cocycle}.} Write
$$
(\Phi)_F=\Phi_D+\hbar^{n+1}\xi+\hbar^{n+2}\psi
\thickspace\mod\hbar^{n+3}
$$
for some $\psi\in\Ug^{\otimes 3}$. Since $\Phi_D$ is
equal to 1 mod $\hbar^2$ and satisfies the pentagon
equation with respect to $\Delta_F(\cdot)=F\Delta
(\cdot)F^{-1}$, we have, mod $\hbar^{n+3}$,
\begin{equation*}
\begin{split}
0
&=
\Pent_{\Delta_F}((\Phi)_F)\\
&=
\Pent_{\Delta_F}(\Phi_D)+
\hbar^{n+1}\dH^{\Delta_F}(\xi)+
\hbar^{n+2}\dH^{\Delta_F}(\psi)\\
&=
\hbar^{n+1}\dH^{\Delta_F}(\xi)+
\hbar^{n+2}\dH\psi
\end{split}
\end{equation*}
where, for any $\eta\in\Ug^{\otimes 3}$,
$$
\dH^{\Delta_F}\eta=
1\otimes\eta-
\Delta_F\otimes\id^{\otimes 2}(\eta)+
\id\otimes\Delta_F\otimes\id(\eta)-
\id^{\otimes 2}\otimes\Delta_F(\eta)+
\eta\otimes 1
$$
is equal to $\dH$ mod $\hbar$. Applying $\Alt_{4}$ to
both sides, and using lemmas \ref{le:alt alt} and \ref
{le:alt sch}, we find, with $\wt\xi=\Alt_3\xi\in
(\bigwedge^3\g)^{\lD}$
\begin{equation*}
\begin{split}
0
&=
\Alt_{4}(\dH^{\Delta_F}\xi)\\
&=
\frac{\hbar}{2}
\sum_{1\leq a<b\leq 4}(-1)^{a+b}
[f^{ab},\wt{\xi}^{ab,1,\ldots,\wh{a},\ldots,\wh{b},\ldots,4}]\\
&=
\frac{\hbar}{2}
\sum_{1\leq a<b\leq 4}(-1)^{a+b}
[f^{ab},
\wt{\xi}^{a,1,\ldots,\wh{a},\ldots,\wh{b},\ldots,4}+
\wt{\xi}^{b,1,\ldots,\wh{a},\ldots,\wh{b},\ldots,4}]\\
&=
-\hbar\schl f,\wt{\xi}\schr
\end{split}
\end{equation*}

where we used the fact that, for any $x\in\Ug$
$$\Delta_F(x)=\Delta(x)+\hbar[f,\Delta(x)]
\thickspace\mod\hbar^2$$
so that
$$\Delta_F(x)-\Delta_F(x)^{21}=2\hbar[f,\Delta(x)]$$
This proves our claim since
$$
\ol{\pi}_D^3\wt\xi=
\Alt_3\ol{\pi}_D^3\xi=0
$$
\halmos

\section{Uniqueness of twists}
\label{ss:coh uniqueness}

Let 
\begin{align*}
\Phi
&=
1+\hbar^2\varphi+\cdots
\in 1+\hbar^2(\Ug^{\otimes 3}\fml)^\g\\
\Phi_D
&=
1+\hbar^2\varphi_D+\cdots
\in 1+\hbar^2(\UgD^{\otimes 3}\fml)^{\gD}
\end{align*}
be two solutions of the pentagon equation \eqref{eq:pentagon}
which are non--degenerate in the sense of definition \ref
{de:non deg}. Contrary to \S \ref{ss:coh existence}, we
do not assume in this subsection that $\Phi_D=\ol{\pi}_
{D}^3(\Phi)$ but merely that
\begin{equation}\label{eq:weak}
\wt{\varphi}_D=\ol{\pi}_D^3(\wt{\varphi})
\end{equation}
where
$$
\wt\varphi=\Alt_3\varphi\in
(\bigwedge^3\g)^{\g}
\quad\text{and}\quad
\wt{\varphi}_D=\Alt_3\varphi_D\in
(\bigwedge^3\gD)^{\gD}
$$
This implies in particular that if $(\cdot,\cdot)$, $(\cdot,\cdot)
_D$ are the bilinear forms on $\g,\gD$ corresponding to $\Phi,
\Phi_D$ via \eqref{eq:non deg} respectively, then $(\cdot,\cdot)
_D$ is the restriction of $(\cdot,\cdot)$ to $\gD$. We denote
the corresponding standard solutions of the MCYBE for $\g$
and $\gD$ by $\rg,\rD$. Let now
$$F_i=
1^{\otimes 2}+\hbar f_i+\cdots\in
1+\hbar(\Ug^{\otimes 2}\fml)^{\lD},
\quad i=1,2$$
be two elements such that $(\Phi)_{F_i}=\Phi_D$. Since $d_H
f_i=0$, we have
$$\wt f_i=\Alt_2f_i\in(\bigwedge^2\g)^{\lD}$$

\begin{theorem}\label{th:coh uniqueness}
Let $F_1,F_2$ be as above and assume that
$\wt{f}_i=\rg\mod(\bigwedge^2\lD)^\lD$. Then,
\begin{enumerate}
\item there exist elements
$$
u\in 1+\hbar\Ug\fml^{\lD}
\quad\text{and}\quad
\lambda\in\hbar\bigwedge^2\cD\fml
$$
such that
\begin{equation}\label{eq:gauge}
F_2=
\exp(\lambda)\cdot 
u\otimes u\cdot F_1\cdot\Delta(u)^{-1}
\end{equation}
\item If
\begin{equation}\label{eq:projection}
\ol{\pi}_D^2(F_i)=1^{\otimes 2},\qquad i=1,2
\end{equation}
$u$ may be chosen such that $\ol{\pi}_D^1(u)=1$. 
\item If
\begin{equation}\label{eq:new theta fix}
F_i^\Theta=F_i^{21},\qquad i=1,2
\end{equation}
then $\lambda=0$ and $u$ may be chosen such that $u
^{\Theta}=u$. $u$ is then unique with this property.
\item If $|\Dg\setminus D|\leq 1$, then $\lambda=0$ and
$u$ is unique up to multiplication by $\exp(c)$ for some
$c\in\hbar\cD\fml$.
\end{enumerate}
\end{theorem}
\proof (i)--(ii) Set
$$f=\rg-\rD$$
and write
$$f_i=\dH g_i+f+\nu_i$$
for some $g_i\in\Ug^{\lD}$, where $\nu_i=\pi^2_D(\wt{f}_i)
\in(\bigwedge^2\lD)^{\lD}=\bigwedge^2\cD$. Then, replacing
$F_i$ by
$$
\exp(-\hbar\nu_i)\cdot
(1-\hbar g_i)\otimes(1-\hbar g_i)\cdot
F_i\cdot
\Delta(1-\hbar g_i)^{-1}
$$
we may assume that
\begin{equation}\label{eq:normalised 1 jet}
F_i=1^{\otimes 2}+\hbar f\mod\hbar^2
\end{equation}
Note that if \eqref{eq:projection} holds, then, by corollary
\ref{co:HC Ho}
$$
0=
\ol{\pi}_D^2f_i=
\dH\ol{\pi}_D^1g_i
$$
and, replacing $g_i$ by $g_i-\ol{\pi}_D^2g_i$, we
may assume that $\ol{\pi}_D^1(g_i)=0$. Similarly, if
\eqref{eq:new theta fix} holds, then
$$
\dH g_i^{\Theta}+f^{\Theta}+\nu_i^{\Theta}=
\dH g_i+f^{21}+\nu_i^{21}
$$
Since $f^{\Theta}=-f=f^{21}$, this yields
$$
\nu_i^{\Theta}=-\nu_i
\quad\text{and}\quad
\dH g_i^{\Theta}=\dH g_i
$$
whence $\nu_i=0$ since $\Theta$ acts as multiplication by
$+1$ on $bigwedge^2\cD\subseteq\\bigwedge^2\h$ and,
replacing $g_i$ by $1/2(g_i+g_i^{\Theta})$, we may assume
that $g_i^{\Theta}=g_i$.\\

We wish now to construct two sequences
$$
v_{n}\in\Ug^{\lD}
\qquad\text{and}\qquad
\mu_{n}\in\bigwedge^2\cD
$$
such that, setting
$$
u_{n}=(1+\hbar^{n}v_{n})\cdots(1+\hbar v_1)
\qquad\text{and}\qquad
\lambda_{n}=\hbar^{n}\mu_{n}+\cdots+\hbar\mu_1
$$
one has
\begin{equation}\label{eq:F2=F1}
F_2=
\exp(\lambda_{n})\cdot 
u_{n}\otimes u_{n}\cdot F_1\cdot\Delta(u_{n})^{-1}
\end{equation}
mod $\hbar^{n+1}$. If \eqref{eq:projection} (resp. \eqref
{eq:new theta fix}) holds, we require in addition that
$\ol{\pi}_D^1(v_{n})=0$ (resp. $v_{n}^{\Theta}=v_{n}$
and $\mu_{n}=0$) for all $n$.\\

By \eqref{eq:normalised 1 jet}, we may set $v_1=0=\mu_1$.
Assume therefore $v_k,\mu_k$ constructed for $k=1\ldots n
$ and some $n\geq 1$. Let $F'_1$ be defined by the right--hand
side of \eqref{eq:F2=F1} so that
\begin{equation}\label{eq:difference}
F_2=F'_1+\hbar^{n+1}\eta\quad\mod\hbar^{n+2}
\end{equation}
for some $\eta\in(\Ug^{\otimes 2})^{\lD}$. One readily
checks that $(\Phi)_{F'_1}=\Phi_D$. Substracting
from this the equation $(\Phi)_{F_2}=\Phi_D$ and
computing mod $\hbar^{n+2}$, we find that
$$
\dH\eta=
\eta^{23}+\id\otimes\Delta(\eta)
-\Delta\otimes\id(\eta)-\eta^{12}=0
$$
Moreover, $\ol{\pi}_D^2\eta=0$ (resp. $\eta^{\Theta}
=\eta^{21}$) if \eqref{eq:projection} (resp. \eqref
{eq:new theta fix}) holds. Thus, $\eta=\dH v+\mu$
for some $v\in\Ug^{\lD}$ and $\mu\in(\bigwedge^2\g)
^{\lD}$ such that
\begin{gather*}
\ol{\pi}_D^1v=0,\\
v^{\Theta}=v\quad\text{and}\quad\mu^{\Theta}=-\mu
\end{gather*}
if \eqref{eq:projection}, \eqref{eq:new theta fix}
hold respectively. Set $v_{n+1}=-v$ and
\begin{equation*}
\begin{split}
F_1''
&=
(1+\hbar^{n+1}v_{n+1})^{\otimes 2}\cdot
F_1'\cdot
\Delta(1+\hbar^{n+1}v_{n+1})^{-1}\\
&=
(1+\hbar^{n+1}v_{n+1})^{\otimes 2}\cdot
\exp(\lambda_{n})\cdot 
u_{n}^{\otimes 2}\cdot
F_1\cdot
\Delta(u_{n})^{-1}\cdot
\Delta(1+\hbar^{n+1}v_{n+1})^{-1}\\
&=
\exp(\lambda_{n})\cdot 
\left((1+\hbar^{n+1}v_{n+1})u_{n}\right)^{\otimes 2}\cdot
F_1\cdot
\Delta\left((1+\hbar^{n+1}v_{n+1})u_{n}\right)^{-1}
\end{split}
\end{equation*}
where the last equality stems from the fact that $v_{n+1}$
is invariant under $\lD$. We have
$$
F_2=\exp(-\hbar^{n+1}\mu)F_1''\quad\mod\hbar^{n+2}
$$
so the inductive step may be completed by setting $\mu_{n+1}
=-\mu$ provided we can show that $\mu$ lies in $\bigwedge^
{2}\cD$. To see this, let
$$
\ol{F}_2=
1+\hbar f+\hbar^2f_2+\cdots+\hbar^{n+1}f_{n+1}
$$
be the truncation of $F_2$ mod $\hbar^{n+2}$ and define
$$
\xi=\xi(f;f_2,\ldots,f_{n+1})\in(\Ug^{\otimes 3})^{\lD}
$$
by
$$
1\otimes\ol{F}_2\cdot\id\otimes\Delta(\ol{F}_2)\cdot\Phi-
\Phi_D\cdot\ol{F}_2\otimes 1\cdot\Delta\otimes\id(\ol{F}_2)=
\hbar^{n+2}\xi\quad\mod\hbar^{n+3}
$$
By lemma \ref{le:deviation}, $\dH\xi=0$ and, by corollary
\ref{co:primary extension}, $\Alt_3\xi=0$ since $\ol{F}_
{2}$ extends to a solution mod $\hbar^{n+3}$. Similarly,
if $\ol{F}''_1$ is the truncation of $F''_1$ mod
$\hbar^{n+2}$, the corresponding error $\xi''$ satisfies
$\dH\xi''=0$ and $\Alt_3\xi''=0$. Since $\ol{F}''_1=
\ol{F}_2+\hbar^{n+1}\mu$ mod $\hbar^{n+2}$ and, for
$n\geq 1$
\begin{equation*}
\begin{split}
\xi''-\xi
&=
\xi(f;f_2,\ldots,f_{n+1}+\mu)-
\xi(f;f_2,\ldots,f_{n+1})\\
&=
f^{23}(\mu^{12}+\mu^{13})+\mu^{23}(f^{12}+f^{13})-
f^{12}(\mu^{13}+\mu^{23})-\mu^{12}(f^{13}+f^{23})
\end{split}
\end{equation*}
we find, using lemma \ref{le:alt=sch}, that $\schl f,\mu
\schr=0$. By theorem \ref{th:sch coh}, this implies that
$$\mu=\schl f,x\schr+y$$
where $y\in\bigwedge^2\cD$ and $x\in\g^{\lD}=\cD\subseteq\h$.
Since $f$ is of weight $0$, $\schl f,x\schr=-\ad(x)f=0$ whence
$\mu=y\in\bigwedge^2\cD$.\\

(iii) Let $u\in 1+\hbar(\Ug\fml)^{\lD}$, with $u^{\Theta}=u$, be
such that
\begin{equation}\label{eq:fix and theta}
u\otimes u\cdot F_1\cdot\Delta(u)^{-1}
=
F_1
\end{equation}
We claim that $u=1$. Assume that $u=1$ mod $\hbar^{n}$ for
some $n\geq 1$ and write $u=1+\hbar^{n}u_{n}$ mod $\hbar^
{n+1}$, where $u_{n}\in\Ug^{\lD}$ is fixed by $\Theta$. 
Taking the coefficient of $\hbar^{n+1}$ in \eqref
{eq:fix and theta} we find that $\dH u_{n}=0$. This
implies that $u_{n}$ lies in $\g$ and therefore in $\h$
since it is of weight zero. Since $\Theta$ acts as $-1$
on $\h$ however, $u_{n}=0$ as claimed.\\

(iv) If $|\Dg\setminus D|\leq 1$, then $\bigwedge^2\cD=0$
so that $\lambda=0$. Let now $u\in 1+\hbar(\Ug\fml)^{\lD}$
be such that
\begin{equation}\label{eq:fix}
u\otimes u\cdot F_1\cdot\Delta(u)^{-1}
=
F_1
\end{equation}
and write $u=1+\hbar u_1$ mod $\hbar^2$. Taking the
coefficient of $\hbar$ in \eqref{eq:fix}, we find that
$\dH u_1=0$ so that $u_1\in\g^{\lD}=\cD$. Now let
$u^{(2)}=u\cdot\exp(-\hbar u_1)=1+\hbar^2u_2
\mod \hbar^2$. Repeating the above argument with
$u^{(2)}$, we find that $u_2\in\cD$ and finally
that there exists a sequence $u_{n}\in\cD$, $n\geq 1$
such that
$$u=
\prod_{n\geq 1}\exp(\hbar^{n} u_{n})=
\exp(\sum_{n\geq 1}\hbar^{n}u_{n})$$
\halmos


\begin{thebibliography}{ZZZ}

\bibitem[Di]{Di} J. Dixmier, {\it Alg\`ebres enveloppantes}.
Cahiers Scientifiques, Fasc. XXXVII. Gauthier--Villars Editeur,
1974.

\bibitem[DS]{DS} J. Donin, S. Shnider, {\it Cohomological
construction of quantized universal enveloping algebras},
Trans. Amer. Math. Soc. {\bf 349} (1997), 1611--1632.

\bibitem[Dr1]{Dr1} V. G. Drinfeld, {\it Quantum groups}, Proceedings
of the International Congress of Mathematicians (Berkeley 1986),
798--820, AMS, 1987.

\bibitem[Dr2]{Dr2} V. G. Drinfeld, {\it  Quasi--Hopf algebras},
Leningrad Math. J. {\bf 1} (1990), 1419--1457.

\bibitem[EK]{EK} P. Etingof, D. Kazhdan, {\it Quantization of Lie
bialgebras. I}, Selecta Math. (N.S.) {\bf 2} (1996), 1--41.

\bibitem[Lu]{Lu} G. Lusztig, {\it Introduction to quantum groups}.
Progress in Mathematics, 110. Birkh\"auser, Boston, 1993.

\bibitem[MTL]{MTL} J. J. Millson, V. Toledano Laredo, {\it Casimir
operators and monodromy representations of generalised braid
groups}, Transform. Groups {\bf 10} (2005), 217-254.

\bibitem[TL1]{TL1} V. Toledano Laredo, {\it A Kohno--Drinfeld theorem
for quantum Weyl groups}, Duke Math. J. {\bf 112} (2002), 421--451.

\bibitem[TL2]{TL2} V. Toledano Laredo, {\it Flat connections and
quantum groups}, The 2000 Twente Conference on Lie Groups. 
Acta Appl. Math. {\bf 73} (2002), 155--173.

\bibitem[TL3]{TL3} V. Toledano Laredo, {\it Quasi--Coxeter
algebras, Dynkin diagram cohomology and quantum Weyl
groups}, {\sf math.QA/0506529}.

\bibitem[TL4]{TL4} V. Toledano Laredo, {\it Quasi--Coxeter
quasitriangular quasibialgebras and the Casimir connection},
in preparation.
 
\end{thebibliography}
\end{document}